\documentclass[11pt]{amsart}
\usepackage{amsmath,amsfonts}
\usepackage{hyperref}
\usepackage{wrapfig}
\usepackage[dvips]{graphicx}
\usepackage[T1]{fontenc}
\usepackage{mathptmx}
\usepackage{color}
\hypersetup{colorlinks=false}
\newtheorem{remark}{Remark}

\newcommand{\ds}{\mathrm{d}s}
\newcommand{\dt}{\mathrm{d}t}
\newcommand{\real}{{\mathbb R}}

\newcommand{\err}{\text{err}}

\begin{document}

  \title{B\'ezier curves that are close to elastica}
  \author[D. Brander]{David Brander}
	\address{Department of Applied Mathematics and Computer Science\\ Matematiktorvet, Building 303 B\\
Technical University of Denmark\\
DK-2800 Kgs. Lyngby\\ Denmark}
\email{dbra@dtu.dk}

  \author[J.A. B{\ae}rentzen]{Jakob Andreas B\ae{}rentzen}
	\email{janba@dtu.dk}
  \author[A-S. Fisker]{Ann-Sofie Fisker}\email{ansofi@dtu.dk}
  \author[J. Gravesen]{Jens Gravesen}\email{jgra@dtu.dk}
  
  \begin{abstract}
We study the problem of  identifying those cubic B\'ezier curves that are close in the $L^2$ norm
to planar elastic curves.   
The problem arises in design situations where the manufacturing process produces elastic curves; these are  difficult to work
  with in a digital environment.  We seek a sub-class of special B\'ezier curves as a proxy.
We identify an easily computable quantity, which we call the $\lambda$-residual
$e_\lambda$, that accurately predicts a small $L^2$ distance.   
 We then identify geometric
criteria on the control polygon that guarantee that a B\'ezier curve
has $\lambda$-residual below 0.4, which 
effectively implies that the curve is
 within $1\%$ of its arc-length
 to an elastic curve in the $L^2$ norm.  Finally we give two projection algorithms that take an input B\'ezier
curve and adjust its length  and shape, whilst keeping the end-points and end-tangent angles fixed, until it is close
to an elastic curve.
  \end{abstract}

\maketitle

\section{Introduction}
\label{sec:Intro}

\begin{figure}[ht]   
\centering
 \begin{tabular}{cc}
   \includegraphics[height=0.1\textheight]{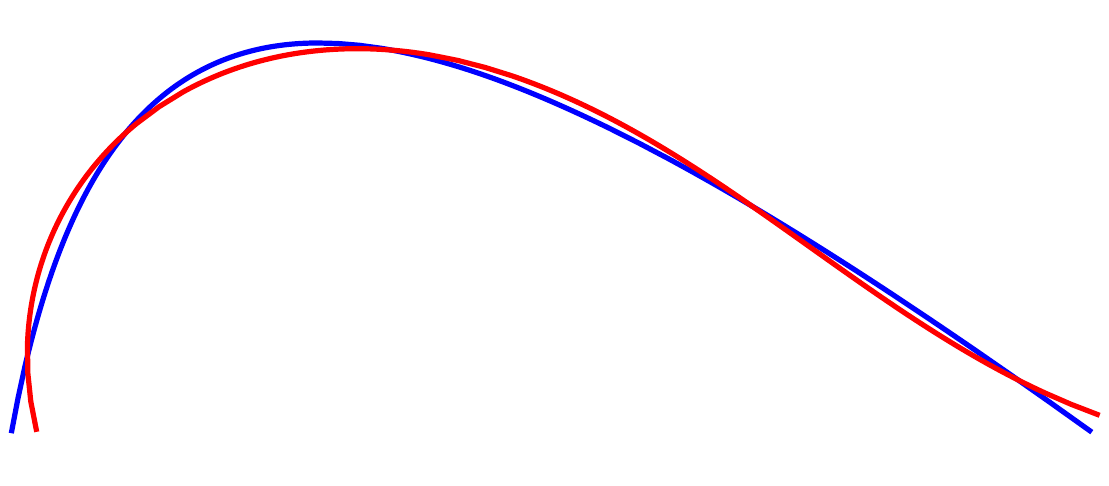} \quad
&  \quad \includegraphics[height=0.1\textheight]{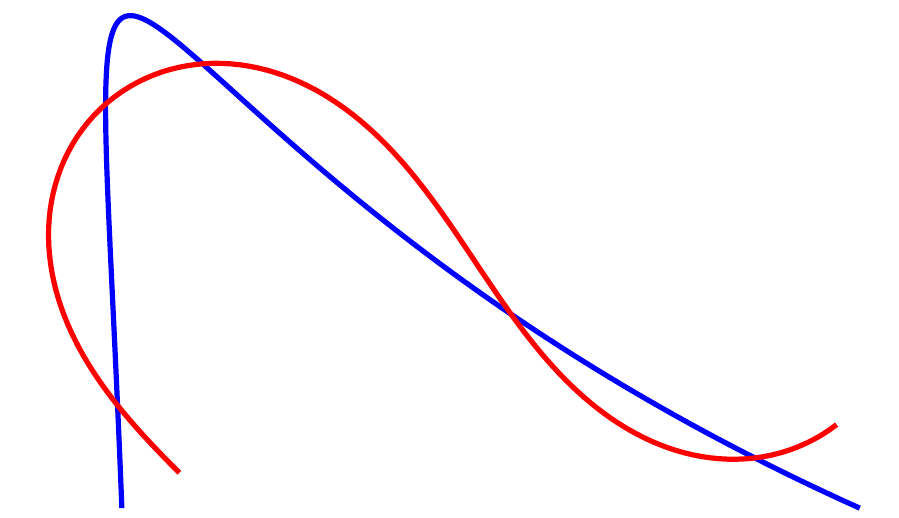}\\
    \footnotesize$\err=   0.006.$ & \footnotesize$\err = 0.032.$ \\
    \end{tabular}    
    \caption{Approximations of two cubic B\'ezier curves (blue) by elastic curves (red). $\err$ is the $L^2$ distance between the B\'ezier and elastic curve.}
    \label{fig:ApproxExamples}
\end{figure}

B\'ezier curves, and their generalization to polynomial and rational splines, were introduced as an easily computable alternative to true splines\footnote{And other templates such as French curves.} around the time that industrial design of ships, aircraft and cars moved  into the digital environment\cite{farin}. True splines, created by thin flexible pieces of wood, held in position at various points, are mathematically described by piecewise planar \emph{elastic curves}, solutions to a nonlinear differential equation, which are difficult to work with compared to polynomials. 
It is sometimes said that a \emph{cubic} B\'ezier curve is quite close to an elastic curve, making them a plausible alternative. This suggestion is based on consideration of the curvature function for a cubic curve that is parameterized by arc-length and not too far from a straight line.  However, it is certainly not true that all cubic B\'ezier curves are close to elastic curve segments:  in Figure~\ref{fig:ApproxExamples}, 
 we have used the algorithm described in \cite{Brander} to obtain approximating elastic curve segments (red) for the given B\'ezier curves (blue).  
Conversely, one can also find elastic curve segments (for example a circle) that are not close to any cubic curve, 
although the space of all B\'ezier curves  is one dimension higher than the space of elastic curve segments (the dimensions are respectively 4 and 3 if scalings, rotations and translations are factored out),
and most shapes produced by elastic splines can be approximated by cubic splines.

From the Computer Aided Design (CAD) point of view, the goal is
not necessarily to replicate exactly the behaviour of true splines, and thus polynomial splines
are usually a good choice. 
However, in some cases there are compelling reasons for faithfully representing a true spline in a digital 
environment: for example if the manufacturing method
naturally produces surfaces swept out by elastic
curves.  An instance of such a method is the recently developed ``hot-blade'' cutting technology \cite{robarch}, whereby architectural formwork is
cut from polystyrene foam using a heated rod, the 
ends of which are controlled by a robot.  \emph{Rationalization} of a CAD design for this production method means segmenting the surface into suitable pieces and then approximating each segment by a family of planar elastic curves.  The primary ingredient for this, an algorithm for approximating an arbitrary curve by an elastic curve segment, is given in \cite{Brander, tokePHD}. 

\begin{figure}
\centering
    \includegraphics[height=0.075\textheight]{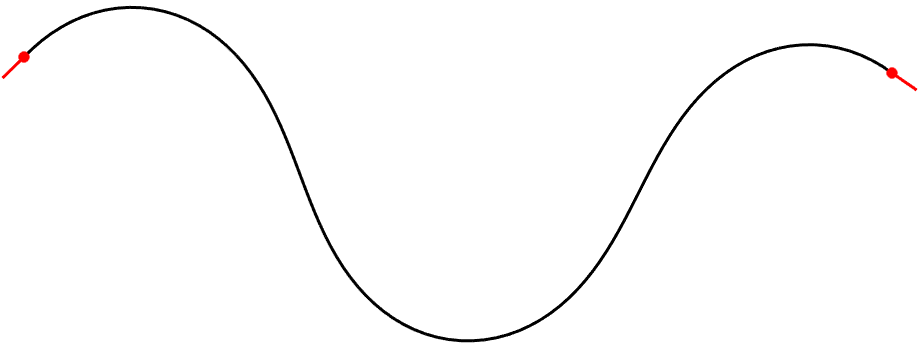} \quad\quad
    \includegraphics[height=0.075\textheight]{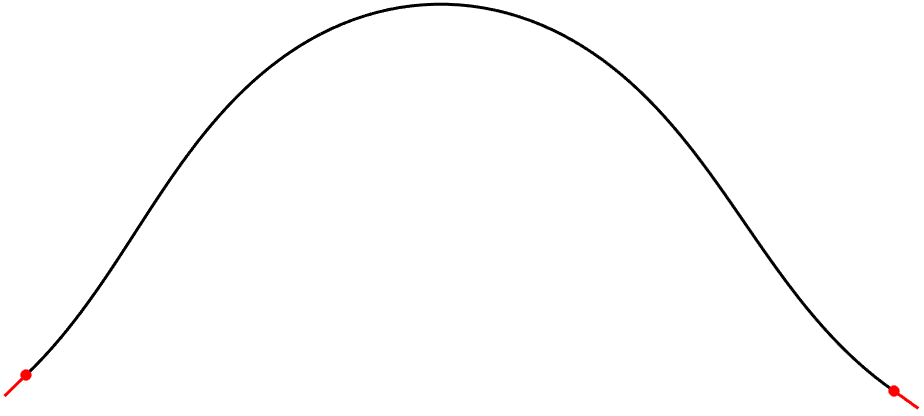}
    \caption{Two elastic curves with identical lengths, end points and tangents.}
    \label{fig:interactive-elastica-challenge}
\end{figure}
As an alternative
to rationalization, we have proposed in \cite{aag2016} a method for using elastic curves themselves in the design process. Unfortunately, true elastic curves can be problematic to work with for a designer. Figure~\ref{fig:interactive-elastica-challenge} shows two elastic curves of the same length,
produced by numerically solving the boundary value problem for elastic curves. Both end points and tangents are (nearly) identical,  but the curves are very different despite both having two inflection points.

One could try to make a tool where the designer is allowed to select between such
 alternatives, or further constraints such as absolute rotation index, number of periods, etc., might be added.
  But  the issue can be side-stepped if a suitable class of cubic B\'ezier curves can be found,
because then the curve is always given uniquely by the control polygon. 
Besides uniqueness there are further benefits to modeling with cubic B\'ezier curves instead of actual elastic curves. Importantly, exchanging data between CAD systems is unavoidable when doing real work, and we posit that modeling with an easily transferable format such as splines is a great benefit.  Along the same lines, an elastic curve would need to be either converted to a spline or a polyline for any downstream usage -- e.g., creating a lofted surface.

In this article, we aim to find conditions on cubic B\'ezier curves such that curves which fulfill these conditions are very close to planar elastic curves. These conditions should be easy to make operational
as algorithms for turning general cubic B\'ezier curves into approximate elastic curves. 
In practical applications,  we can then easily obtain precise numeric solutions for the true elastica as a last step.

\subsection{Related literature}
The energy minimizing (subject to end constraints) property of elastic curves makes them fair from a design point of view, as it is
one of the ways to enforce a smooth curvature function.  This led to various attempts to emulate this property
within a computational setting.  
One idea is to work with some kind of numeric approximation for elastic curves, e.g.\ \cite{Mehlum}, \cite{DBLP:journals/cad/BrunnettK94}, \cite{Edwards:1992:EEN:146847.146925}, \cite{Delas}.
The drawback of this approach is that it is computationally expensive to solve the boundary
value problem associated with an elastic curve; 
furthermore, as mentioned above, there are issues of non-uniqueness causing
instability.  This makes an interactive design tool difficult.

A more practical approach is to replace the elastic curves with another class of curves such as energy minimized quadratic \cite{AHR2014}, cubic \cite{YongCheng}, or quintic \cite{Lu2015} splines, and Pythagorean-hodograph curves \cite{FAROUKI1996227, FaroukiCAGD2016}.
All of these are designed to have low bending energy subject to Hermite interpolation conditions.

Our work does not fit precisely into any of the above viewpoints, because we are not
just interested in working with curves of low bending energy, but rather curves that are (visually) close to \emph{actual}
elastic curve segments.  This is because our practical motivation is to provide a fast, interactive, way to approximate physical elastic curves within a CAD system. Therefore, unlike in previous works, our measure of ``goodness'' for a curve is not its bending energy, but the distance from the closest true elastic curve segment.

Finally, the reader is no doubt aware that,
in geometric design,  there is also a lot of interest in other measures
of fairness, especially minimal curvature variation (\cite{farin2008}),
among others (e.g., \cite{WuYang2016}, \cite{Yanetal}, \cite{MG2014}). 
The driver of these works is purely the aesthetic appearance of the curves, in contrast to our
to work, which is motivated by fabrication constraints, namely the
physical proximity to a curve segment that is an elastic rod.

\subsection{Overview}

In Section~\ref{sec:elastic-curves} we briefly introduce results from \cite{Brander} on approximating
an arbitrary curve by an elastic curve and an important quantity that we call the \emph{$\lambda$-residual}, 
$e_\lambda$. This easily computable quantity measures how close the curvature function of a given curve 
is to the curvature function of an elastic curve.  In Section \ref{sec:experiments} we take a large sample space
of cubic B\'ezier curves and show that the $\lambda$-residual is sufficiently well correlated with the $L^2$-distance 
from the curve to an elastic curve to allow us to use $e_\lambda$ as a proxy for this distance.

The goal is then to find a method of projecting an arbitrary cubic B\'ezier curve to a B\'ezier curve
that has $e_\lambda$ below some threshold.  For the sake of concreteness, we aim for $e_\lambda <0.4$,
which corresponds approximately to an $L^2$ distance less than 0.007, geometrically 
a deviation of 0.7 percent relative to the length of the curve.
In this work we concentrate on projections that alter only the length of the B\'ezier curve, keeping
the endpoints and end tangent angles fixed. This is grounded in the idea that a tool that allows
the user to prescribe exactly the end-points and end-tangents of a curve is most useful for a designer.

In Section \ref{sec:gradient} we briefly discuss the possibility of gradient driven approaches to projection,
i.e., minimizing either the elastic energy or the $\lambda$-residual $e_\lambda$ subject to fixed end-data.
 This approach is
problematic due essentially to the existence of multiple local minima.

In Section \ref{sec:geometric} we find a large subset $\Pi$ of cubic B\'ezier curves, characterized by the
end tangent angles, the inner polygon angles and the lengths of the two outer polygon edges, such
that $e_\lambda$ is bounded by $0.4$ on $\Pi$.   We also describe an end-data preserving algorithm
for projecting into $\Pi$.   The results of this section are included partly because the geometric 
characterization has the possibility of being adapted to other projections if the end-tangent preserving
requirement is dropped.  In Section \ref{resultssubsection} we check the result
by taking a large random sample of curves inside $\Pi$ and finding an $L^2$-approximation of each
curve by an actual elastic curve segment, with the result that 100\% of these B\'ezier curves have
an $L^2$-distance below 1\% of arc-length from an elastic curve segment.

Finally, in Section \ref{sec:Feedback} we present a significantly more effective end-data preserving
projection algorithm based on 
computing the value of $e_\lambda$ in real time.

\section{Elastic curves and the $\lambda$-residual}
\label{sec:elastic-curves}

\begin{figure}[ht]
    \centering
    \includegraphics[width=0.82\textwidth]{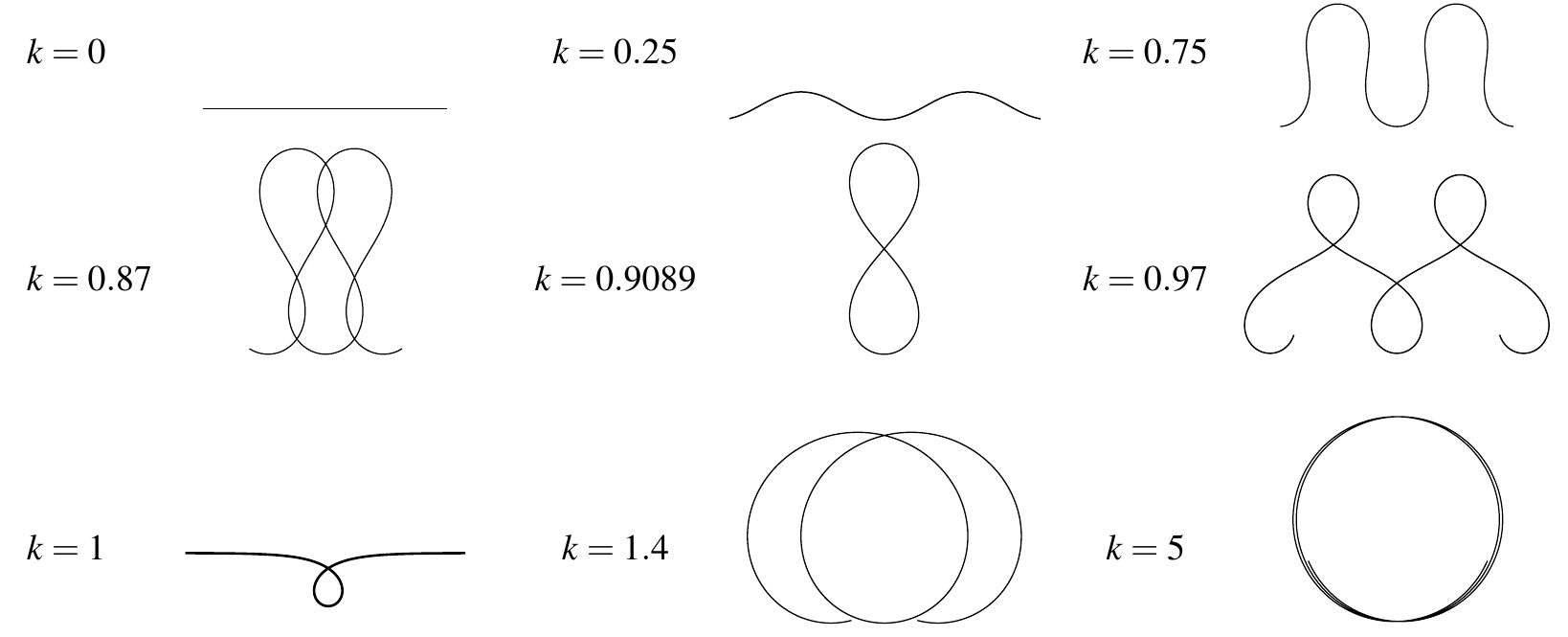}
    \caption{Examples of elastic curves.}
    \label{fig:elcrv}
\end{figure}

An elastic curve is defined to be the minimizer of the \emph{elastic energy}
$\int \kappa^2(s) \,\ds$ among curves with given
endpoints, end-tangents and length.  We refer the reader to 
\cite{Brander} for the relevant theory, as well as for the 
details of the approximation algorithm mentioned here.
This algorithm uses the well-known analytic description of an elastic curve: 
\begin{equation*}
\label{elas}
\xi_{k}(s) = (2E(s,k)-s,2k\,(1-\text{cn}(s,k)))\,, 
\end{equation*}
where $\text{cn}$ is a Jacobi elliptic  function,\footnote{
Note that this representation uses
an extension of the elliptic functions and integrals to $k \in [1,\infty)$
(see \cite{Brander}), which allows 
a single formula incorporating the
elastic curves both with and without inflections. The limit as $k \to \infty$ is a circle.}
$E$ is the incomplete elliptic integral of the second kind and $k\in[0,\infty)$.
Then any elastic curve segment is, up to
a scaling, rotation and translation, given by a piece of some $\xi_{k}$.
Some examples of these curves are shown in Figure~\ref{fig:elcrv}.

The approximation algorithm described in \cite{Brander} takes as input any curve and returns an approximating elastic curve segment. The algorithm consists of two steps:
\begin{enumerate}
   \item Obtain a \emph{first guess} elastic curve that has roughly the same shape as the input curve. This exploits
	the fact that the curvature of an elastic curve is affine in a particular direction, and the first guess is obtained by
	solving linear least squares problems. 
    \item  Apply an optimization to adjust the parameters of the first guess to obtain an optimal solution.
\end{enumerate}
We use this algorithm in this article to obtain 
an approximating elastic curve segment $\gamma_e$ for a given B\'ezier curve $\gamma_B$. In the computation of the first guess we find the parameters $\lambda_{1},\lambda_{2},\alpha\in \mathbb{R}$  that  minimize $\int^{1}_{0} (\kappa(t) +\lambda_{1}y(t) -\lambda_{2}x(t)-\alpha)^{2}(\ds/\dt)\,\dt$ (this is equivalent to solving a linear system). The associated residual, \emph{$\lambda$-residual},  is given by the formula:
\begin{equation*}
    \label{lambdares}
    e_{\lambda} = \sqrt{\int^{1}_{0} (\kappa(t) +\lambda_{1}y(t)-\lambda_{2}x(t)-\alpha)^{2}\frac{\ds}{\dt}\,\dt}
    \Biggm/
    \sqrt{\int_{0}^{1}\kappa(t)^{2} \frac{\ds}{\dt}\,\dt}\,,
\end{equation*}
where $\kappa$ is the curvature of the B\'ezier curve, $s(t)$ is the arclength of the B\'ezier curve at $t$.

The $\lambda$-residual   can be taken as a measure of how much the curvature of a given curve
$(x,y)$ deviates from being the curvature of an elastic curve, and we will use this as one of our tools to analyze B\'ezier curves.  
Because the $\lambda$-residual can be
computed at interactive speeds, it is the most important tool we have for measuring how
 close a curve is  to an elastic curve.

\section{ The $\lambda$-residual as measure of closeness between cubic B\'ezier curves and elastic curves }
\label{sec:experiments}
\begin{wrapfigure}{r}{0.5\textwidth}
  \begin{center}
  \includegraphics[height=60mm]{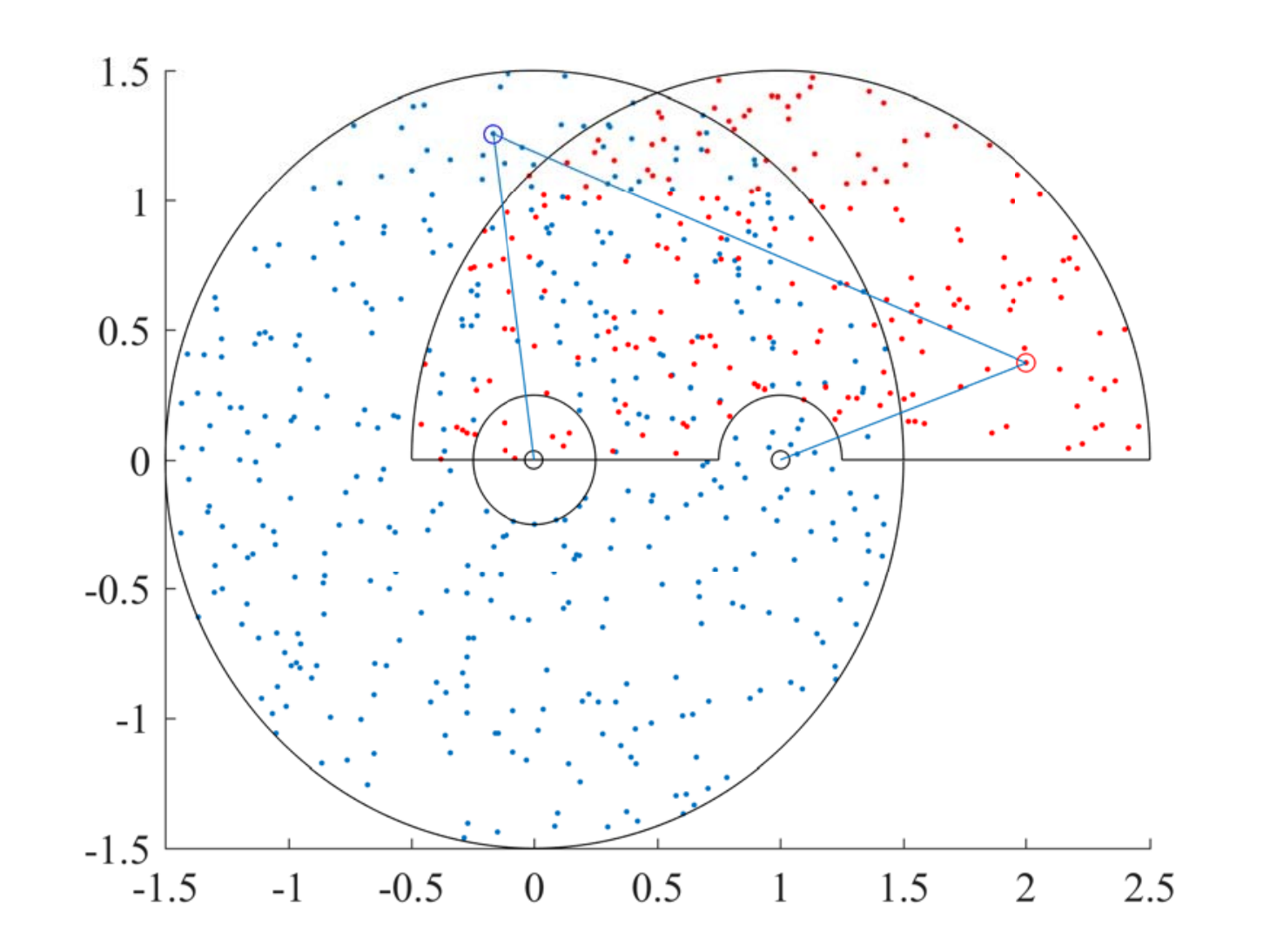}
	\end{center}
  \caption{Initial sample space of inner control point pairs: one of the 80,000 control polygons 
	is shown for illustration.}
  \label{fig:bez}
	\end{wrapfigure}
As our starting point we study a randomized 
sample of 80,000 cubic B\'ezier curves to determine usable criteria 
for identifying those that are close to elastic curves. 
We take as our initial yardstick the $L^2$ distance (see \eqref{L2dist} below) 
to the approximating elastic curve, 
which is obtained by the approximation algorithm in \cite{Brander}.
 We begin with the criterion that the good B\'ezier curves are those with  an $L^2$ distance less than $0.01$ (i.e., $1$ percent of the curve length) from the approximating elastic curve. To ensure the quality of the optimal solution we have applied the elastic curve approximating algorithm using different optimization routines in MATLAB: trust-region, interior point, SQP and Quasi-Newton.   \\

\begin{figure}[ht]
   \centering
   \unitlength = 50mm
      \includegraphics[height=\unitlength]{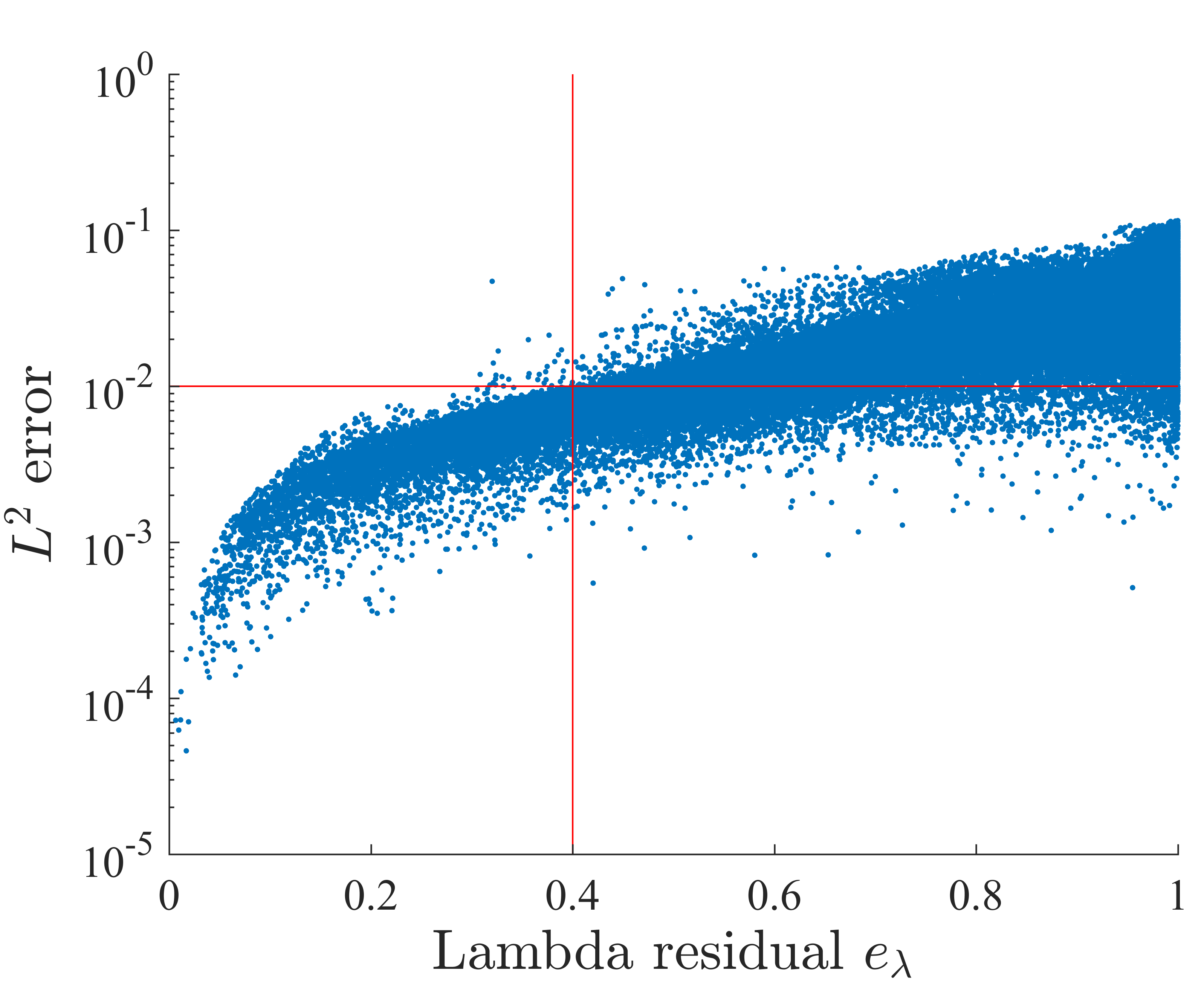} \quad
			\includegraphics[height=\unitlength]{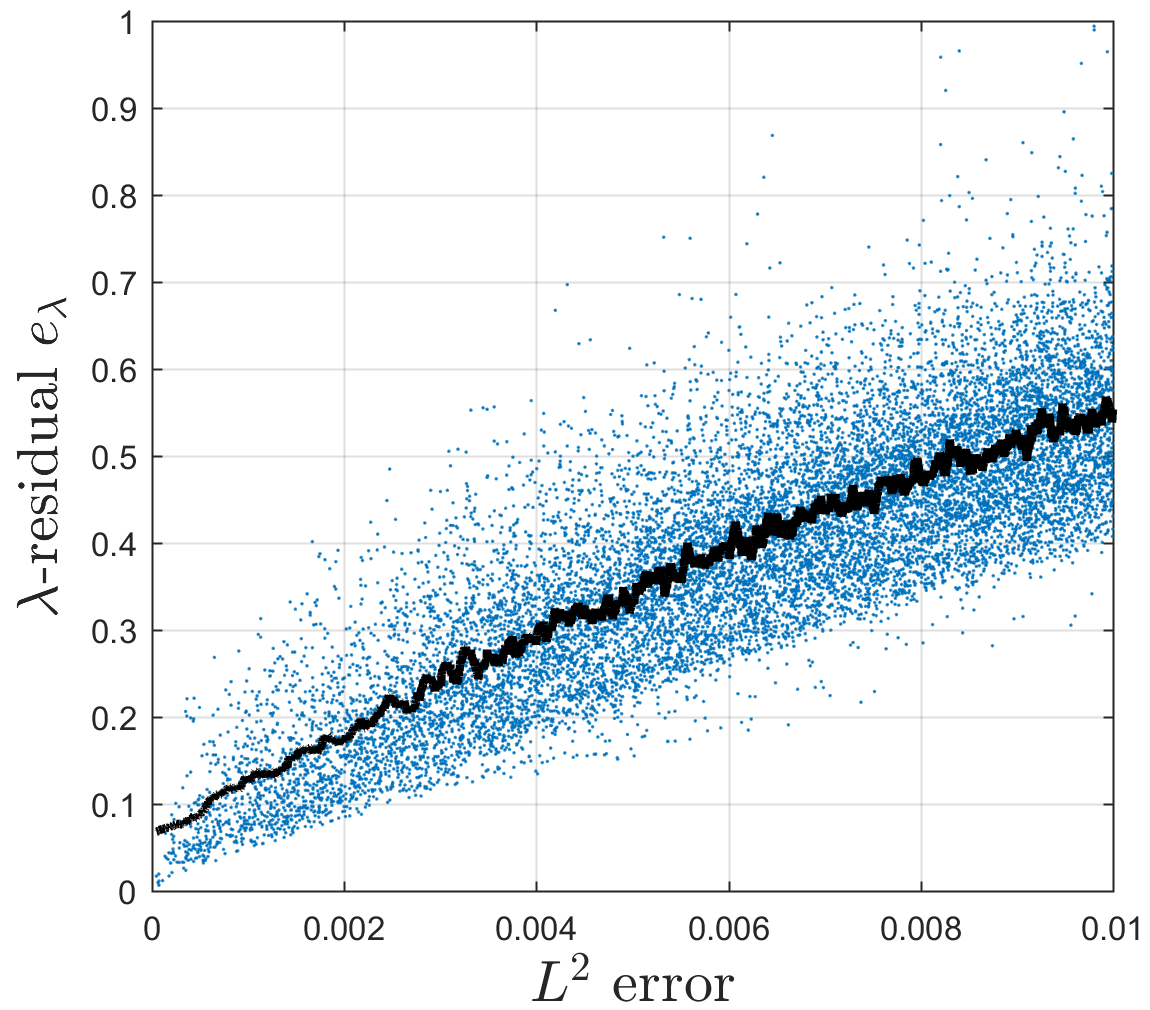} 
   \caption{The $\lambda$-residual and $L^2$ error.}
   \label{fig:scatplot4}
\end{figure}

All of our arguments are invariant under scaling, rotation and translation, so we can fix the endpoints of the B\'ezier curve at $(0,0)$ and $(1,0)$ and only vary the two inner control points. Our initial sample space consists of all possible B\'ezier curves obtained from pairing the red and blue points  in  Figure~\ref{fig:bez} and using them as the inner control point
pairs.
We found approximating elastic curves for all the curves in the  sample space,
and computed the $L^2$ error.  We then plotted the $L^2$ distance
against  the $\lambda$-residual, $e_\lambda$,
defined in Section \ref{sec:elastic-curves}.  From Figure~\ref{fig:scatplot4} we conclude that, for B\'ezier curves, the $\lambda$-residual is well correlated with the $L^2$ distance to an elastic curve, especially for small values.
A maximum $L^2$ error of $0.01$ is achieved by almost all of the curves with $e_\lambda<0.4$.

Compared to finding an approximating elastic curve and its  $L^2$ error (via optimization), 
the $\lambda-$residual is easy to compute:  we only have to solve a linear system to obtain $\lambda_1,\lambda_2$ and $\alpha$. Furthermore the $\lambda$-residual depends continuously on the control polygon vertices of a spline curve, and is more reliable because the $L^2$ distance is subject to such factors as the optimization used
to find the approximating curve. We therefore use the $\lambda$-residual as a
more practical measure of closeness to elastic curves. From the trend curve in 
Figure \ref{fig:scatplot4}, we apply the following measure of quality:
\begin{table}[ht]
    \centering
    \begin{tabular}{ l r r r }
   &         Best   &   Good   & Borderline  \\
$L^2$ error maximum  &  0.003      & 0.007     &  0.01
\\
$\lambda$-residual   maximum  & 0.22 & 0.4 & 0.5 
\\ 
\end{tabular}
\end{table}

\section{Projection of cubic B\'ezier curves:  problems with the gradient driven approach}
\label{sec:gradient}

The $\lambda$-residual  can tell us if a cubic B\'ezier curve is close to an elastica. This can be used as a diagnostic tool in a design framework. However, if we want to model elastic curves with cubic B\'ezier curves we need a projection tool for the curves that have a large $\lambda$-residual. 

In this article we consider only projections that keep the end-points and end-tangents of the
curve fixed: that is, given an input B\'ezier curve, we will modify it only by moving the two
inner control points along the line segments between them and their corresponding end-points.

\begin{figure}[ht]
   \centering
   \unitlength = 50mm	   
   \includegraphics[height=\unitlength]{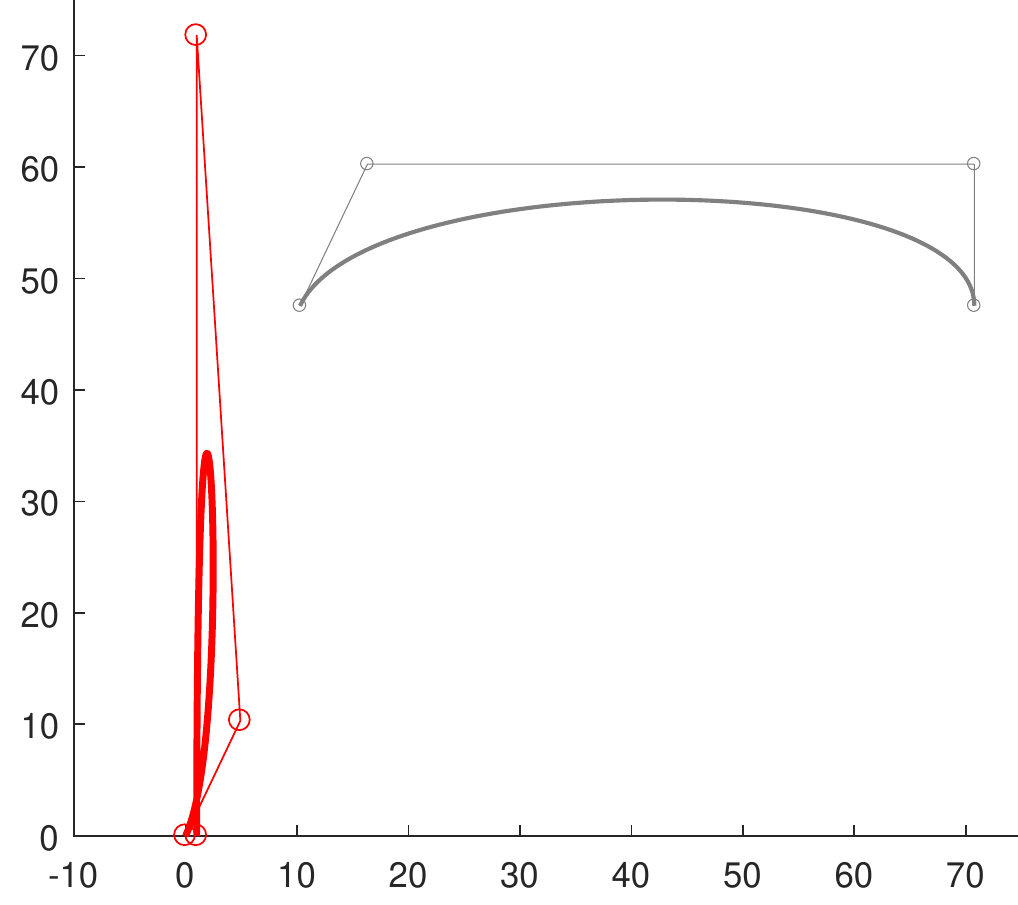} 
	      \quad \quad \quad
\includegraphics[height=\unitlength]{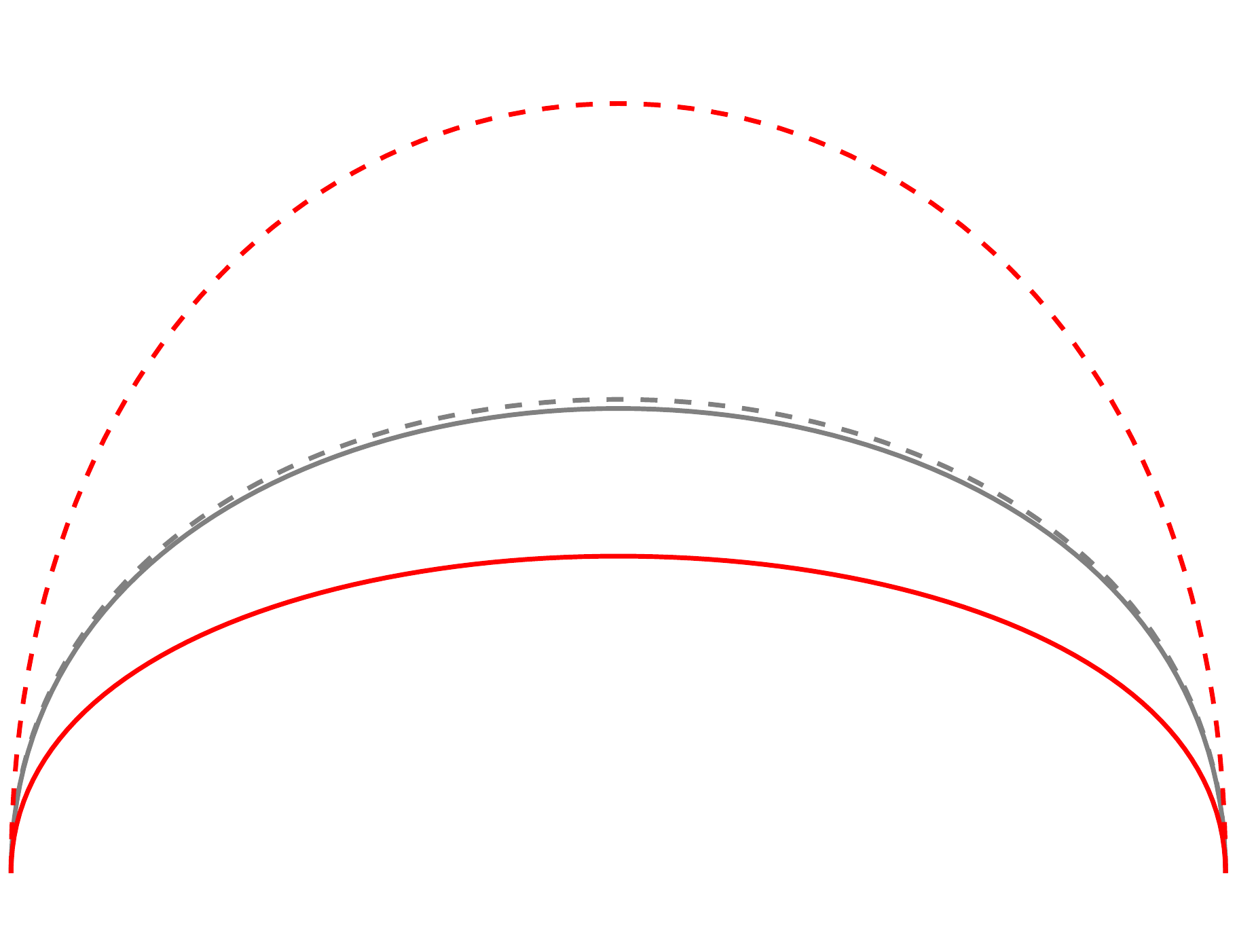} 
   \caption{Left: input (grey), output (red, different scale) after some steps of minimizing the bending energy with fixed end tangent
	directions. Right: two similar inputs (grey) and very different output (red) after minimizing the $\lambda$-residual. }
   \label{fig:gradexample}
\end{figure}

We first consider a gradient driven projection.  
We have two candidates for the energy to minimize:
\begin{enumerate}
\item The bending energy $\int \kappa^{2} \ds$.
\item The $\lambda$-residual $e_{\lambda}$.
\end{enumerate} 
We have implemented and tested a projection tool for both energies.  In both cases,
we keep the end-points and end-tangent directions fixed and allow the length to vary.
With the bending energy, except in the case of a straight line, 
the length of the curve will keep increasing, as this reduces the elastic energy, see Figure~\ref{fig:gradexample} left.  This means that the projection does not return a useful B\'ezier curve.  
For both energies we have an issue that there is more than one local minimum
for the optimization:  for example, for given end-points and end-tangent angles, there are often
local minimizers to be found both with and without polygon intersections, as well as 
minimizers corresponding to both inflectional and non-inflectional elastic curves.  Figure~\ref{fig:gradexample}, right, shows an instance of this, where the $\lambda$-residual energy was used.  From these observations we conclude that the gradient driven approach is  unreliable, and/or does not always return a B\'ezier curve that is close to an elastic curve.

\section{A projection based on a geometric characterization of curves with low $\lambda$-residual}  \label{sec:geometric}
Due to the problems with the gradient driven approach we now look for a projection tool based on a more geometric  characterization of cubic B\'ezier curves that are close to elastic curves. 
We will first find criteria in terms of the edge lengths
and polygon angles, that guarantee that  $e_\lambda \leq 0.4$.

\subsection{Finding the projection zone}
\begin{figure}[ht]
   \centering
   \unitlength = 40mm
	      \includegraphics[height=\unitlength]{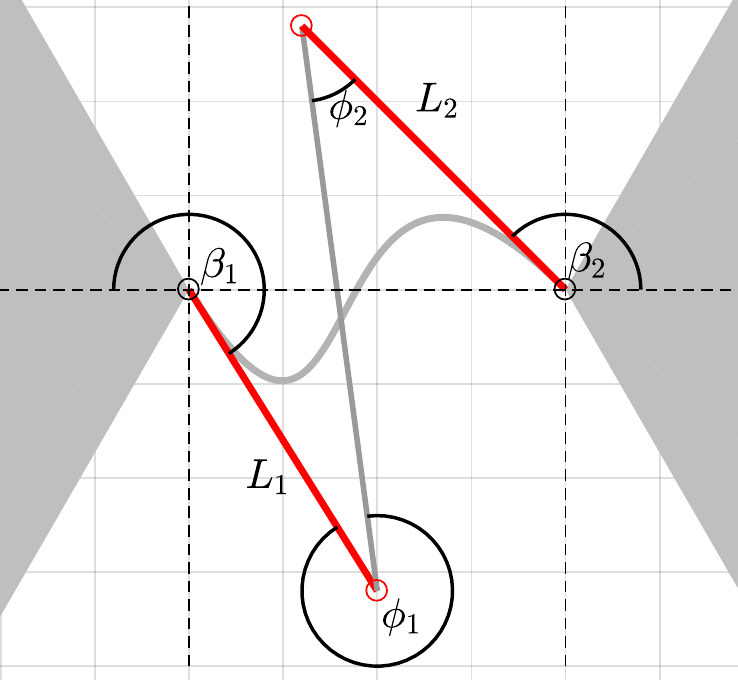} \quad \quad 
	   \includegraphics[height=\unitlength]{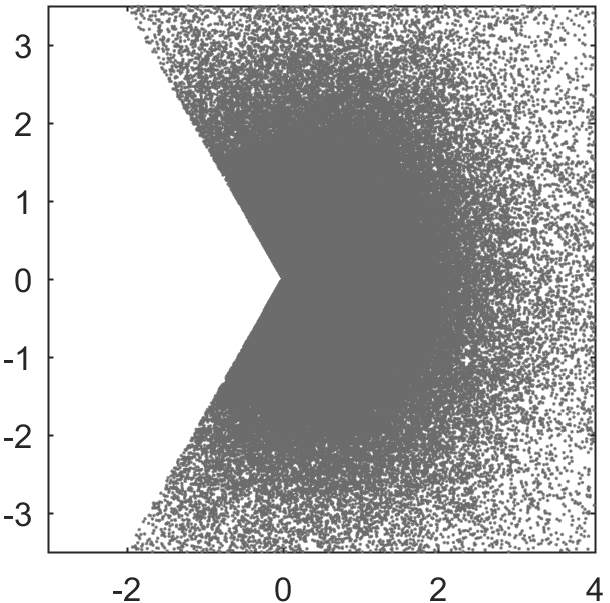} \quad
		   \includegraphics[height=\unitlength]{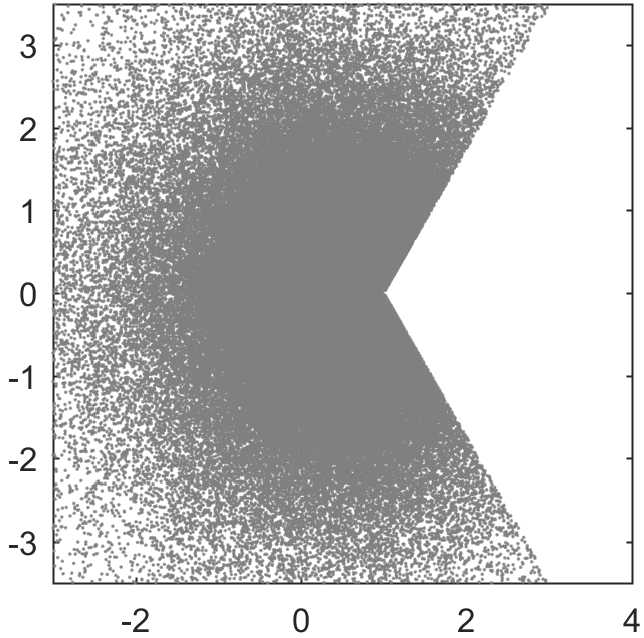}
   \caption{Left: Angle constraints: the shaded region is excluded. Extremely asymmetric shapes are
	excluded by a bound on $|\beta_1-\beta_2|$. Right: the sample space consists of B\'ezier curves with endpoints at
	$(0,0)$ and $(1,0)$ and middle control points all possible pairs with $p_1$ chosen from the first set shown,
	and $p_2$ chosen from the second. }
   \label{fig:angle_constraints}
\end{figure}

We first made  a large, randomly distributed sample of quadruples of points, $p_0$, $p_1$, $p_2$ and $p_3$ in the unit disc. Then we scaled and rotated each polygon with $p_0 \mapsto (0,0)$ and $p_3 \mapsto (1,0)$, to get
a random collection of B\'ezier curves in standard position. We reduced the sample space by removing: curves with self-intersecting polygons, and curves that do not satisfy the 
following angle constraints (see Figure~\ref{fig:angle_constraints}): \\
\[
\hbox{ \textbf{Angle Constraint 1 (absolute angle constraint):}} \quad
 \beta_1, \, \beta_2 \in \left(  \pi/3 \, , \,  2\pi -\pi/3 \right),
\]
where $\beta_1 \in [0,2\pi]$ is the angle measured clockwise from the negative $x$-axis to
the first polygon edge, and $\beta_2$ is the symmetric analogue 
 for the third edge measured anticlockwise from the positive $x$-axis.\\
\[
\hbox{\textbf{ Angle Constraint 2 (symmetry constraint):} } \quad| \beta_1 - \beta_2| < 0.4\pi.
\]
 The angle constraints are applied because there are relatively few B\'ezier curves with low $\lambda$-residual that do not satisfy them. These constraints were arrived 
 at by statistical analysis; but one can see how they arise by considering the shape of B\'ezier curves that do not satisfy them. A curve that fails 
 Constraint 1 will have very high curvature in some region
unless  the curve is extremely
 symmetric.  Curves that fail Angle Constraint 2 are inflectional, and the more the constraint
is exceeded, the higher will be the curvature variation.
With the two angle constraints, we are left with 4,330,509 curves containing  434,580 of the
``best'' curves with $e_\lambda<0.22$, but 3,187,526 with $e_\lambda>0.4$, these latter being 
either unacceptable or borderline unacceptable curves.

\begin{figure}[htb]
    \centering
    \begin{tabular}{cc}
        \includegraphics[height=60mm]{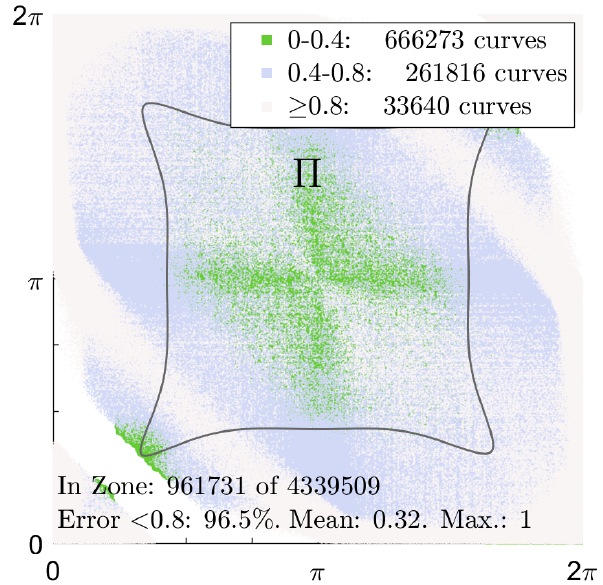} \quad & \quad
            \includegraphics[height=60mm]{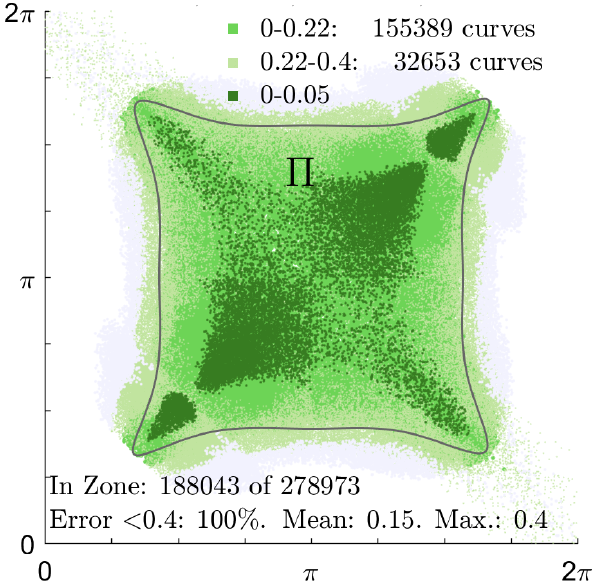}  \\
						No length constraints & With length constraints  
        \end{tabular}
    \caption{Scatter plot of $(\phi_1, \phi_2)$ values for the curves in the angle-constrained sample space.
		Numbers refer to the curves \emph{inside} the projection zone $\Pi$ bounded by the dark blue curve. }
    \label{fig:projzone1}
\end{figure}

Let $\phi_1, \phi_2 \in [0, 2\pi)$ be the angles from the first and third control polygon edges to the second,
as shown in Figure~\ref{fig:angle_constraints}, left, and $L_1$ and $L_2$ denote the lengths of the
first and third edges respectively.
We look for conditions on these quantities that will remove the bad curves, whilst keeping a large number of  good curves.   Figure \ref{fig:projzone1} (left) shows a plot
of $(\phi_1, \phi_2)$ from the sample space, colored by the $\lambda$-residual.
This shows that $(\phi_1,\phi_2)$ alone cannot be used to characterize the good curves. We next
 include lower and upper bounds on $L_1$ and $L_2$: \\
\[
\hbox{\textbf{ Edge-Length Constraint 1:}}  \quad L_{min} \leq L_1, \,\, L_2 \leq L_{max},
\]
where
\begin{eqnarray*}
L_{min}=\hbox{max}(0.4(1+6\Delta), \,\, 0.27), \\
   L_{max}=\hbox{max}(1.2(1+5\Delta),  \,\, 0.58),  \\
		\Delta= \hbox{sign}(\theta_2)\frac{(\theta_1-\theta_2)}{\pi},
	\end{eqnarray*}
	$\theta_i$ are the two end-tangent angles, measured from the positive $x$-axis,
and we limit the relative lengths of the outer polygon edges: \\
\[
\hbox{\textbf{  Edge-Length Constraint 2 (Symmetry for inflectional curves):}} \quad
\hbox{max} \left (L_1/L_2, \, L_2/L_1\right) \leq 1.3,
\]
for curve  $\phi_2 \leq \pi \leq \phi_1$ or $\phi_1 \leq \pi \leq \phi_2$. 

The edge-length constraints were chosen in an ad-hoc manner, based on the observation that we need $L_1$ and $L_2$ to be shorter if the legs are angled inward toward the center of the polygon, and longer if angled outwards, and that, for inflectional curves, the better curves have outer edges of similar length. \\

With the edge-length constraints we obtain a bounded region $\Pi$ (see Figure \ref{fig:projzone1}, right)
 that \emph{only} contains good curves, with a maximum $\lambda$-residual of $0.4$, whilst retaining $36\%$ of the best curves ($e_\lambda<0.22$). We  use this region as a geometrically delineated
set of cubic B\'ezier curves that are close to elastic curves. We define the boundary of $\Pi$ to be
 the closed curve given by 
taking the points $(1.05, 1.05)$, $(1.9, 1.3)$, $(\pi, 1.35)$, $(4.3,1.3)$,  $(5.2,2\pi-5.2)$,
reflecting them about the line $\phi_1=\phi_2$, and then about the line $\phi_1=2\pi-\phi_2$,
and interpolating using a periodic cubic spline interpolation.\\

\subsection{End-tangent angle preserving projection to $\Pi$}
\begin{figure}[htb]
    \centering
    \begin{tabular}{cc}
        \includegraphics[height=40mm]{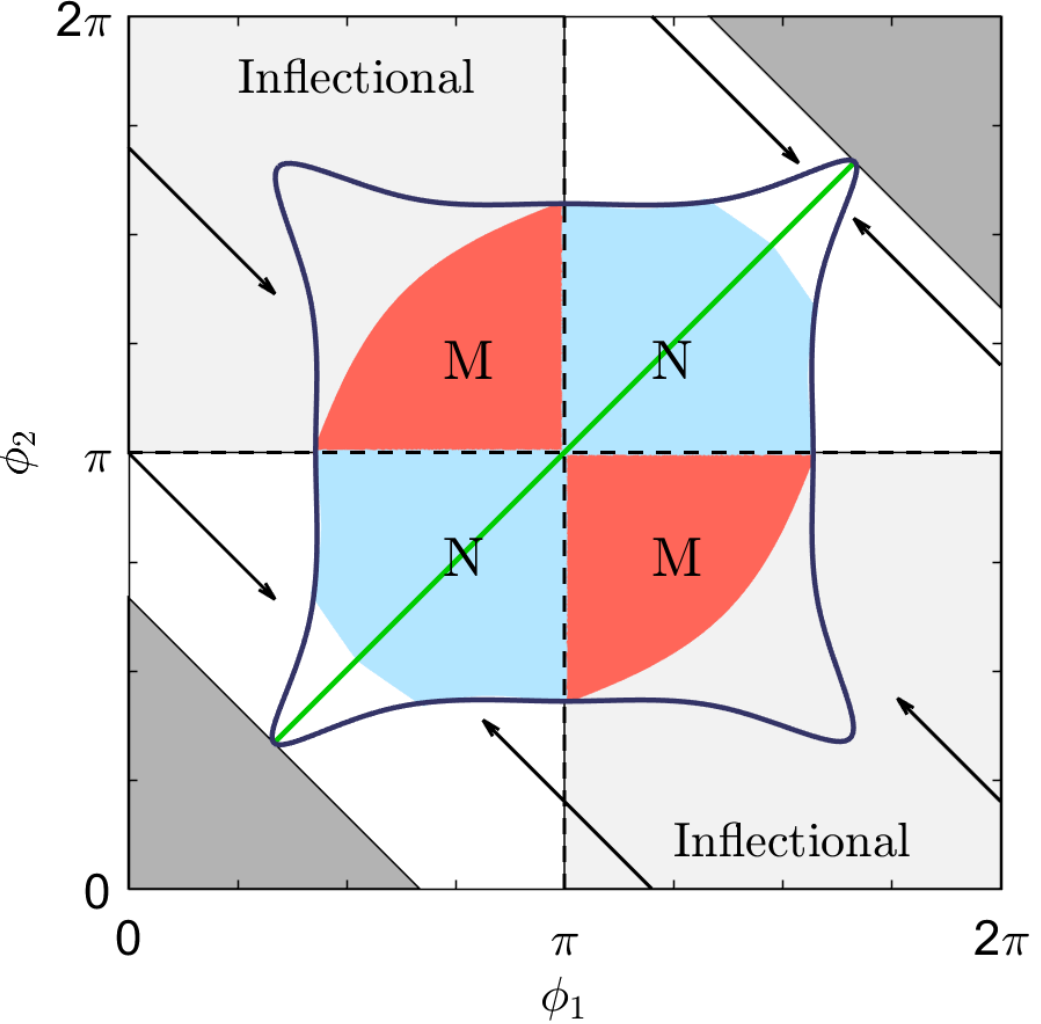}  \quad &
				   \includegraphics[height=38mm]{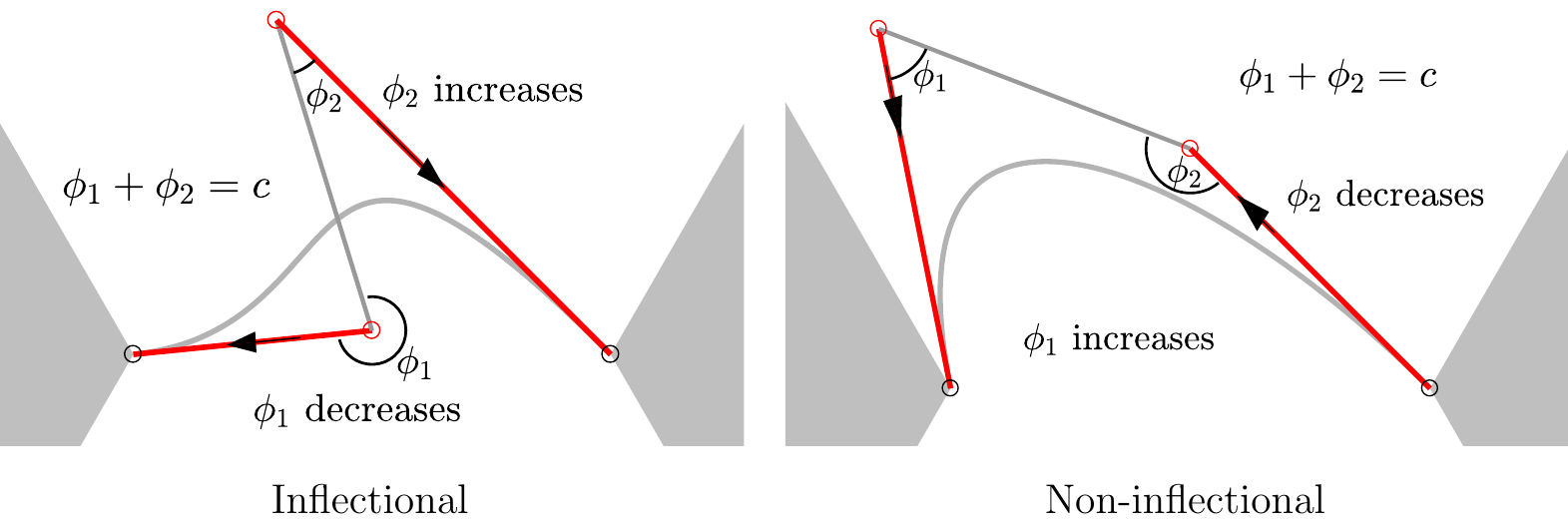} 
        \end{tabular}
    \caption{End-tangent angle preserving projection. The sum $\phi_1+\phi_2$ remains constant.
		For non-inflectional curves, the vertex with the larger angle is moved away from its end-point.}
    \label{fig:projection}
\end{figure}

We now describe, for an input B\'ezier curve that satisfies Angle Constraints 1 and 2,
a projection onto $\Pi$ that preserves both endpoints and end tangent angles.
The only adjustments made to the input curve is to move the
two inner control points along the lines between themselves and their corresponding end points
(see Figure \ref{fig:projection}). 
It follows from elementary geometry that the effect of such a motion on the pair $(\phi_1, \phi_2)$  of inner polygon angles is to move the pair up and down the line $\phi_1+\phi_2 = \hbox{constant}$ on which it lies. We choose to move always in the direction of the arrows in the diagram to the left in Figure \ref{fig:projection}, i.e., towards the line $\phi_1=\phi_2$, and it follows that a curve with an inflection can change to a curve without an inflection, but not the other way around.

The projection algorithm  is as follows: first remove any self-intersection of the control polygon edges by
shortening the outer edge-lengths $L_1$ and $L_2$.  Next adjust $L_1$ and $L_2$ if necessary so 
that the edge-length constraints are satisfied (this will not introduce self-intersections). The B\'ezier curve
is then inflectional or non-inflectional depending on the values of $(\phi_1, \phi_2)$, and the 
following routines are applied for the two cases:\\
\noindent \textbf{ Inflectional Curves ($\phi_2 \leq \pi \leq \phi_1$ or $\phi_1 \leq \pi \leq \phi_2$):}
We make changes that \emph{decrease} either $L_1$ or $L_2$ or both using
the deformation 
\begin{equation}\label{infproj}
L_1(t)=(1-t)L_1+tL_{\min},  \quad L_2(t) = (1-t)L_2 + t L_{\min},
\end{equation}
with $t \in [0,1].$
This moves  $(\phi_1,\phi_2)$ towards the line $\phi_1=\phi_2$. 
The minimal edge length constraint is $ L_{\min} = 0.27$ for all relevant inflectional curves,
 and the locus of $(\phi_1,\phi_2)$ that 
satisfy $L_1=L_2=L_{\min}$ is plotted as the red region $M$  in Figure \ref{fig:projection}, left.
This shows that we are guaranteed either to reach the projection zone $\Pi$, or that
the curve becomes non-inflectional, \emph{before} we arrive at the minimal edge-length constraint.\\

\begin{figure}[htb]
    \centering
    \begin{tabular}{ccc}
        \includegraphics[height=32mm]{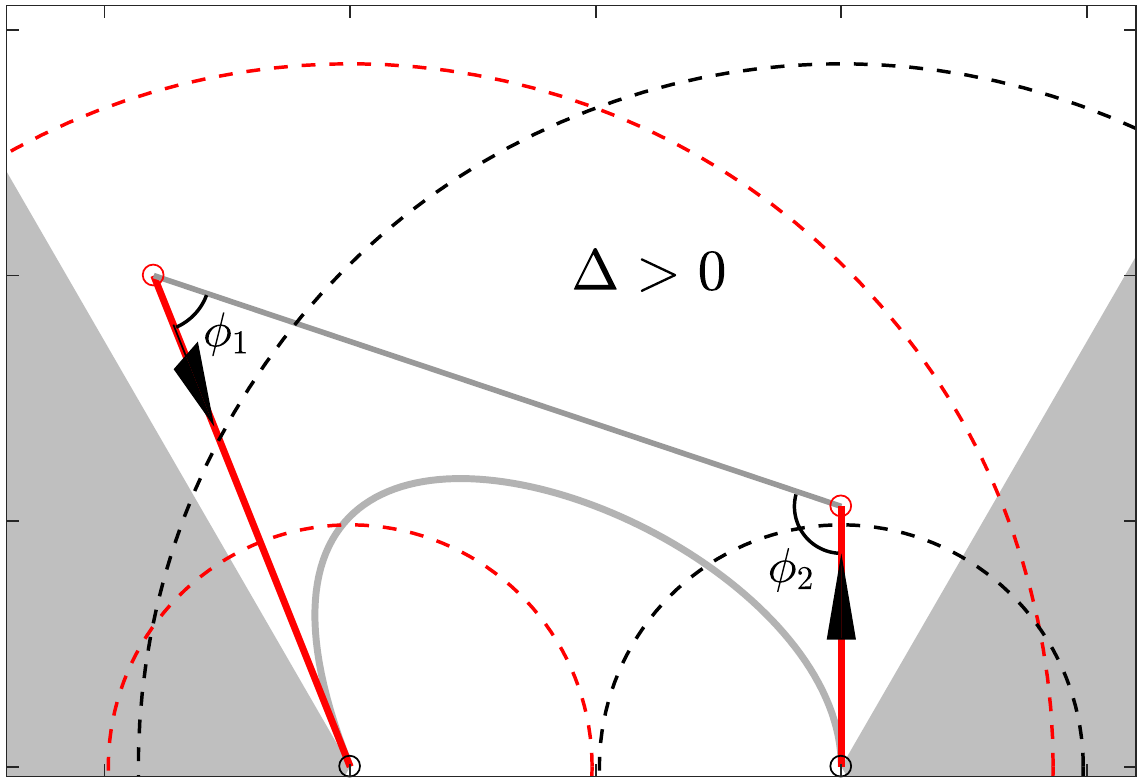}  \quad &
				   \includegraphics[height=32mm]{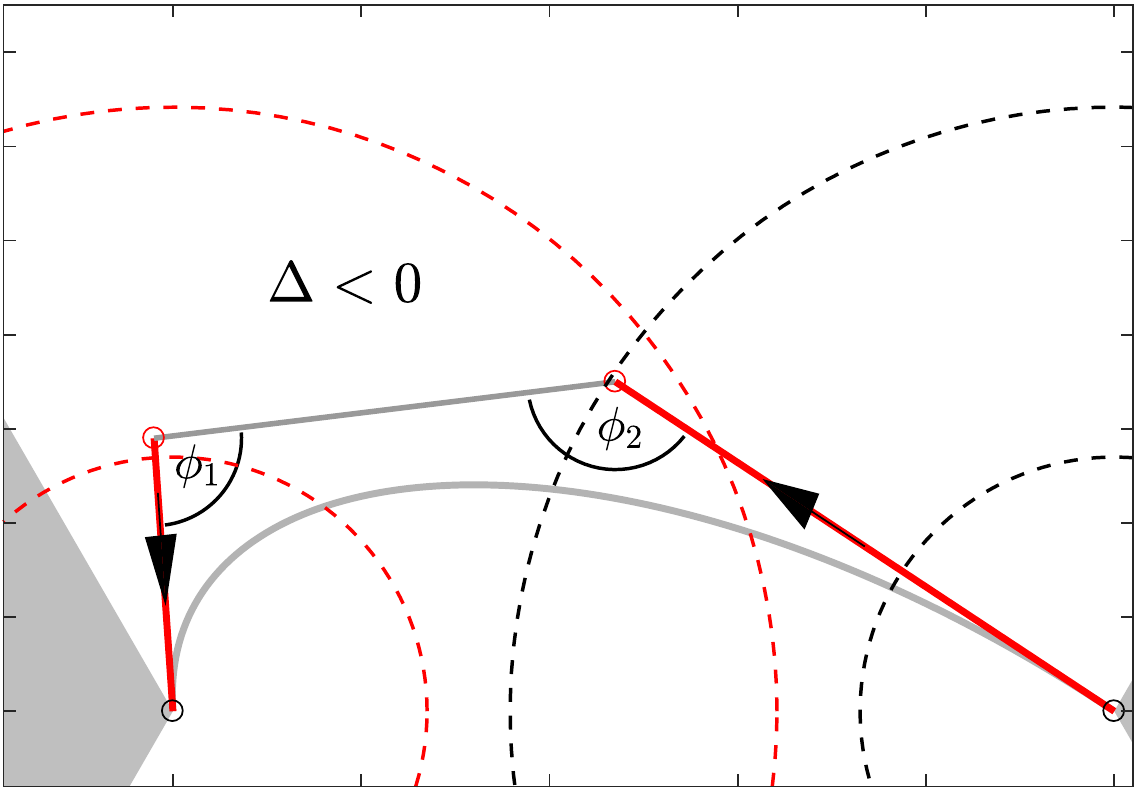}  \quad &
					   \includegraphics[height=32mm]{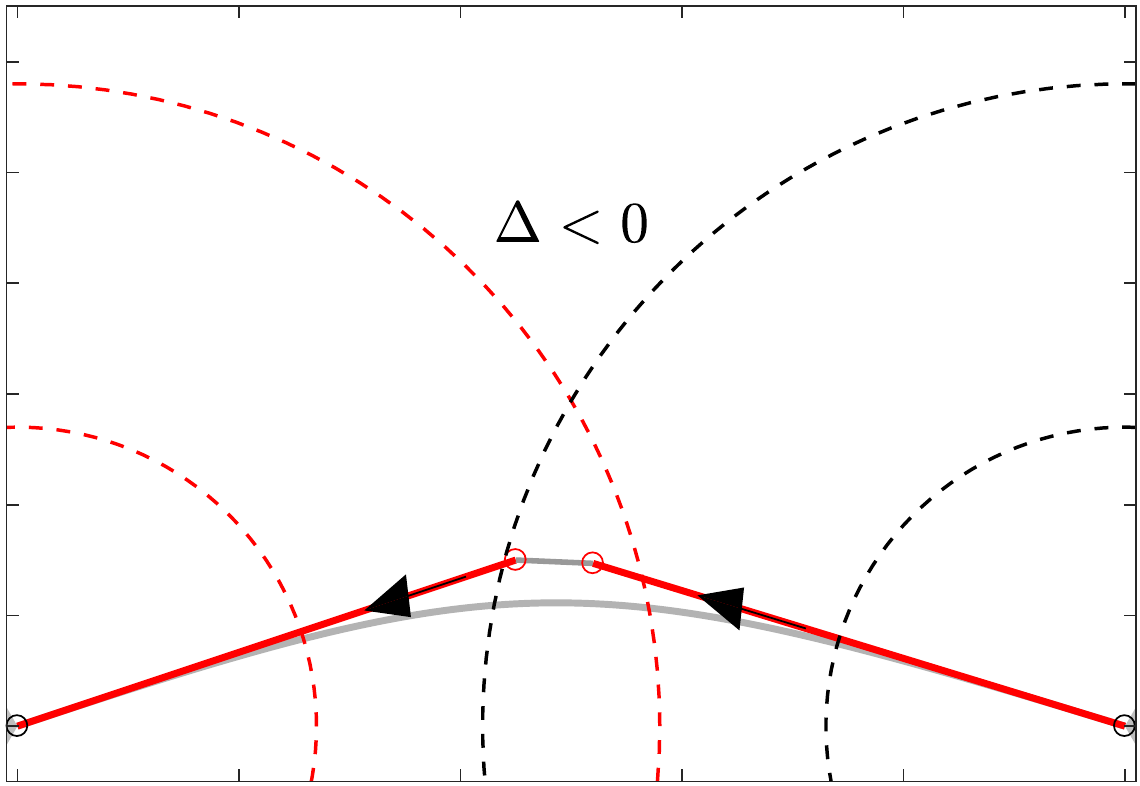} 
        \end{tabular}
    \caption{Length constraints on the two outer edges when projecting non-inflectional curves.}
    \label{fig:projnoninfl}
\end{figure}

\noindent \textbf{Non-inflectional Curves ($\phi_1, \, \phi_2 \leq \pi$ or $\phi_1, \,\phi_2 \geq \pi$):}\\
Adjust the angles $(\phi_1, \phi_2)$: we describe the case $\phi_1<\phi_2 \leq \pi$, (the other cases are
analogous). Either decreasing $L_1$ or
increasing $L_2$ (or both) will move $(\phi_1, \phi_2)$ towards the line $\phi_1=\phi_2$
 (Figure \ref{fig:projection}, right image).
We decrease $L_1$ and increase $L_2$ simultaneously, e.g., with the formulae
\begin{equation} \label{noninfproj}
L_1(t)=(1-t)L_1+tL_{\min},  \quad L_2(t) = (1-t)L_2 + t L_{\max},
\end{equation}
with $t \in [0,1]$, until the projection zone is reached.    
As in the inflectional case, we need to be sure that the projection zone is reached before the
edge-length constraints are violated (See Figure \ref{fig:projnoninfl}). 
The projection could only fail if $\phi_1<\phi_2$ and if $L_1=L_{\min}$
and $L_2 = L_{\max}$ (or the analogue for the other positions).  To verify that this does not happen
 we calculate the locus of all such points, and this is plotted as the blue region $N$
in Figure \ref{fig:projection}.  Since $N$ is contained inside $\Pi$, the edge-length
limits are not exceeded here either.

 Finally, note that it is not difficult to verify that no self-intersections are introduced by the above procedures.   \\

\subsection{Results of the purely geometry based approach.}  \label{resultssubsection}
The above projection algorithm is easily implemented by discretizing the parameter $t \in [0,1]$,
and  looping until the projection zone is reached.
The projection zone  is a compact subset of $\real^4$,
consisting of a bounded subset of the pairs of middle control points.  The $\lambda$-residual
depends continuously on the control points: we computed the $\lambda$-residual for a new set
of $10.7$ million densely distributed points in the projection zone, and found a maximum of $0.393$, confirming the results in Figure
\ref{fig:projzone1}.   Thus we are, with a high level of confidence,
 guaranteed a result below $0.4$ with this method.   

\begin{figure}[htb]
    \centering
    \begin{tabular}{ccc}
        \includegraphics[height=40mm]{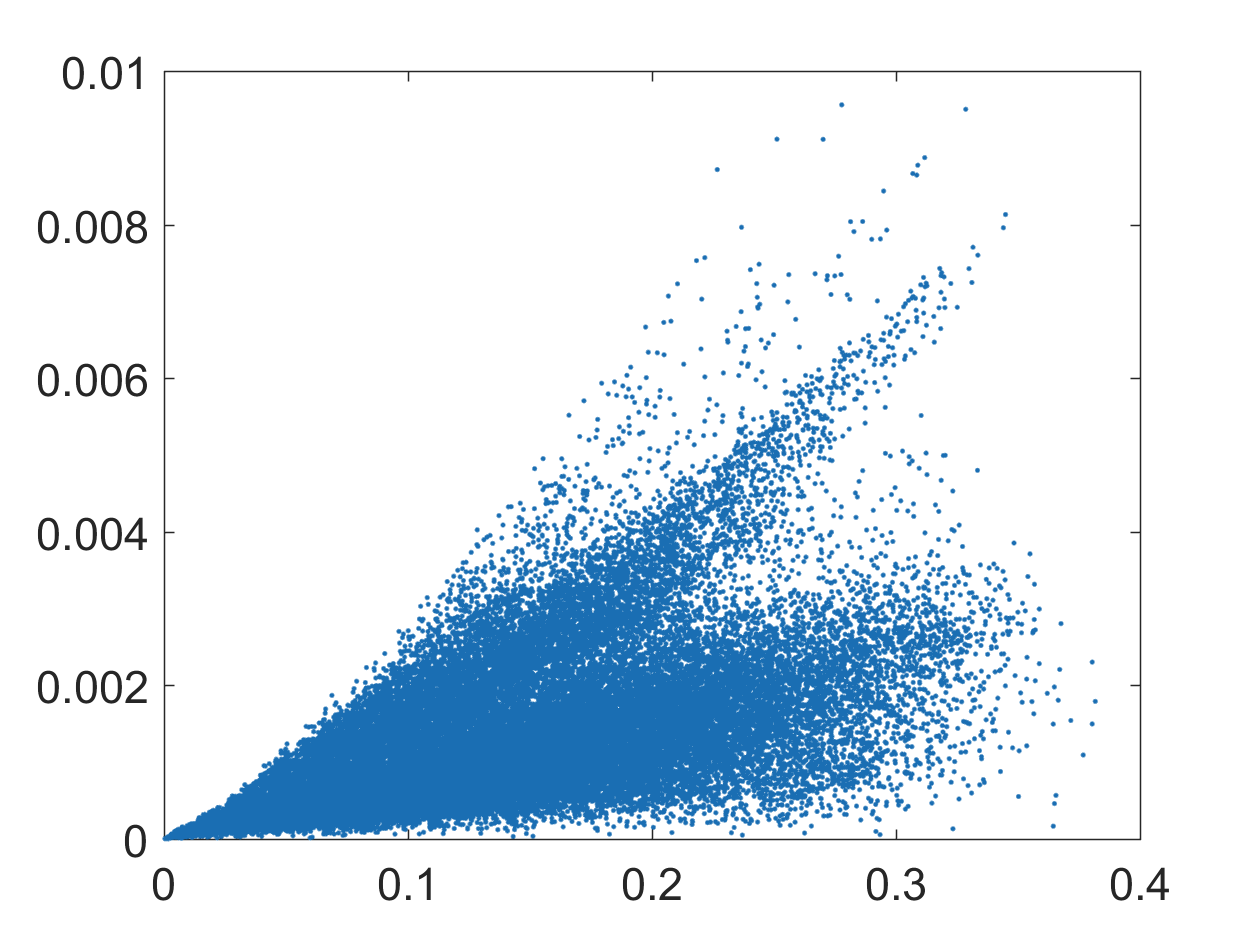}  & \quad
				   \includegraphics[height=40mm]{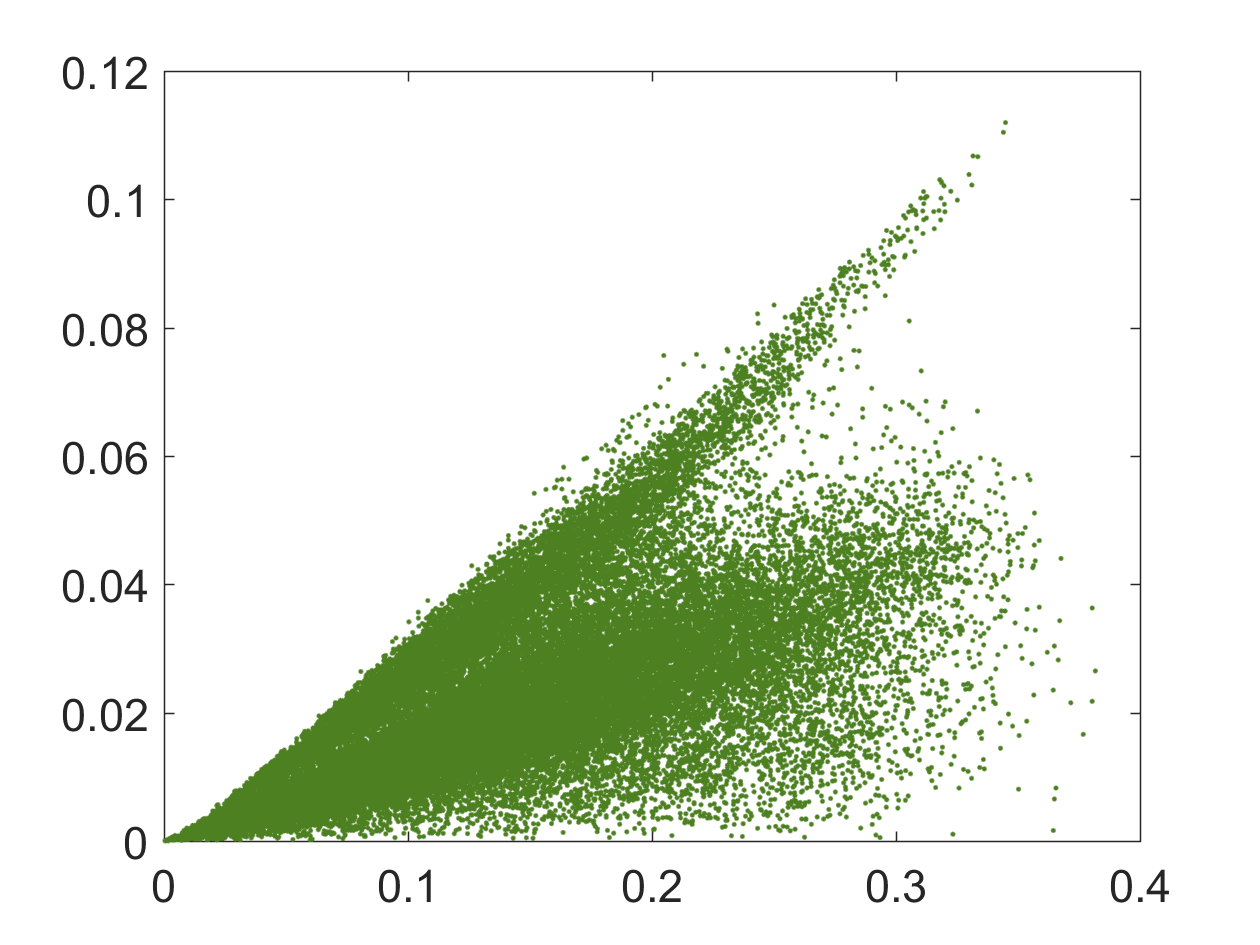} &
        \end{tabular}
    \caption{Results of optimized $L^2$ approximations (left) and $H^1$ approximations (right) for curves in the projection zone.  In each case the distance is plotted against $e_\lambda$ of the input curve.} 
    \label{fig:L2H1}
\end{figure}

Finally, to confirm the quality of the curves in the projection zone, we took a random sample of 38,826 curves in this
zone and approximated them using the method of \cite{Brander},  both with the $L^2$ and the
$H^1$ distances between the B\'ezier curve $\gamma_B$ and the elastic curve $\gamma_e$:
\begin{equation} \label{L2dist}
\textbf{$L^2$ distance : } \sqrt{\int_0^{1} \frac{\|\gamma_B(t)-\gamma_e(s(t)/L)\|^{2}}{L^{3}}\, \|\gamma_B^\prime(t)\| \, \dt} \, ,
\end{equation}
\begin{equation*}
\textbf{$H^1$ distance : } \sqrt{\int_0^{1} \frac{\|\gamma_B(t)-\gamma_e(s(t)/L)\|^{2}}{L^{3}}\,
\|\gamma_B^\prime(t)\| \, \dt + \int_0^{1} (\theta_B(t)-\theta_e(s(t)/L))^{2}\, \frac{\|\gamma_B^\prime(t)\|}{L}  \, \dt} \, ,
\end{equation*}
where $\theta_B$ and $\theta_e$ denote respectively the tangent angles of $\gamma_B$ and $\gamma_e$, 
$L$ is the length of the curve $\gamma_B$, and $s(t)$ is the arc-length function for  the curve $\gamma_B$.  Note that theses norms are invariant under scaling, and are relative to the length
of the B\'ezier curve.
The results, plotted against $e_\lambda$, are shown in Figure \ref{fig:L2H1}.    For the $L^2$ distance, they show that the result is as expected, i.e., the maximum $L^2$ distance is always less than 0.01, and almost always below 0.007.
The $H^1$ distance shows that an upper bound on $e_\lambda$ can also be used to obtain an
upper bound on the $H^1$ distance from an elastic curve;  this means that we could also use $e_\lambda$
as a proxy for the $H^1$ distance.

\section{A projection algorithm based on feedback}
\label{sec:Feedback}
Using the projection described above, we obtain a  B\'ezier curve with a maximum value $e_\lambda = 0.4$.
However,  one finds that  $98\%$ of curves score below $0.3$. Moreover, there are also curves with 
very low values of $e_\lambda$ \emph{outside} the projection zone.    
Since the $\lambda$-residual can be computed easily at interactive speeds, we can obviously improve both the quality of
the result (lower value of $e_\lambda$) as well as the size of the design space (no need to move all the
way into $\Pi$ in many cases), by computing $e_\lambda$ on the fly and applying the deformations
\eqref{infproj} and \eqref{noninfproj}
until  either a threshold level for $e_\lambda$ is reached or $t=1$.

Analysis shows that, if we have fixed the end-tangent angles but allow the lengths $L_1$ and $L_2$ to vary,
then the minimal value of $e_\lambda$  is close to the boundary of region $M$ for the inflectional curves, and hence it makes sense to use \eqref{infproj} for the inflectional case.  
For non-inflectional curves the minimal value of $e_\lambda$ does not occur exactly along the line $\phi_1-\phi_2=0$,
so it turns out that we can do better than using \eqref{noninfproj} for these.   

\subsubsection*{Feedback-based projection algorithm}
\begin{enumerate}
\item Choose a threshold level $E \in [0,1)$ for $e_\lambda$ and a minimum allowed length $L_{min}$ for the outer edges of inflectional curves.  We obtained the  best results by choosing $L_{min}=2.7$.
\item Scale and rotate so that the control points $p_0$, $p_1$, $p_2$ and $p_3$ are in standard position, 
with end-points $p_0=(0,0)$ and $p_3=(1,0)$. Denote this transformation by $T$.
\item Remove any self-intersection of the polygon edges by reducing the lengths $L_1$ and 
$L_2$.
\end{enumerate}
Now we classify the curve and adjust $L_1$ and $L_2$:
\begin{enumerate}
\item[4a.]
The input  curve is classified as \textbf{inflectional} if we obtain an inflectional curve after adjusting the outer edges to
be within the range $[L_{min}, L_{max}]$, where $L_{max}$ is defined by the same formula as in Edge-Length Constraint 1.
 In this case, we take the adjusted curve as the input curve, then
run the deformation formula \eqref{infproj} until either (a) both $L_1$ and $L_2$
are equal to $L_{min}$,  or (b) $e_\lambda \leq E$, or (c) the curve becomes non-inflectional.
\item[4b.]
For \textbf{non-inflectional curves}, we choose suitable terminal values $\mathcal{L}_1$ for 
$L_1$ and $\mathcal{L}_2$ for $L_2$ (see below).  Then adjust $L_1$ and $L_2$ according to the formula
\begin{equation} \label{noninfproj2}
L_1(t)=(1-t)L_1+t \mathcal{L}_1,  \quad L_2(t) = (1-t)L_2 + t  \mathcal{L}_2
\end{equation}
with $t \in [0,1]$, until either (a) $L_1 = \mathcal{L}_1  $ and $L_2 =\mathcal{L}_2$, 
or (b) $e_\lambda \leq E$. 
\item[5.] Finally, apply $T^{-1}$ to the new curve to obtain the projected B\'ezier curve.
\end{enumerate}
The values $\mathcal{L}_i$ are chosen as local minimizers for $e_\lambda$. They depend on the end tangent angles, so should be expressed 
as $\mathcal{L}_i(\theta_1, \theta_2)$, where $\theta_i$ are the angles of the outer edges to the
$x$-axis. 

We chose, for $\theta_1 \geq 0$ (i.e. if $y_1 \geq 0$), 
\begin{align*}
	\mathcal{L}_1(\theta_1, \theta_2)  :=  & \min(0.12, f(\theta_1, \theta_2)), \\
 f(\theta_1, \theta_2) = &0.00001475 \exp(10.39 \theta_1 -10.48 \theta_2) 
  + 0.4574 \theta_1  \exp(1.711 \theta_1 -2.535 \theta_2)  \\
&	+ 2.772 \theta_2 \exp(-0.08504 \theta_1 -0.9109 \theta_2) 
	-0.2957 \theta_1 \theta_2 \exp(-0.6606 \theta_1),\\
	\end{align*}
For curves with $y_1 \leq 0$, the corresponding formula is obtained by symmetry, and the formula for $\mathcal{L}_2$ is also the symmetric analogue. We obtained $f$ by computing local minimizers for a discrete set of $(\theta_1,\theta_2)$ given by non-inflectional curves that satisfy Angle Constraints 1-2. We then fitted the data using this formula giving a small sum of squares error. 

Note: In practice, the algorithm is implemented by discretizing the interval $[0,1]$.  If the threshold 
$E$ is not reached before the other terminating conditions, then the $t$-value corresponding to the lowest
value of $e_\lambda$ is chosen as the solution.  This ensures that any inaccuracy in the choice of the function $f$
has minimal impact on the result.

\begin{remark}
There is no unique or best choice of $\mathcal{L}_i(\theta_1, \theta_2)$, except in the case that the
input curve is unambiguously close to either an inflectional or non-inflectional elastic curve.  Many input curves are close to
both types of elastica, and this leads to two different local minima for $e_\lambda$.
Therefore, we have made some arbitrary choices in our definition of $f$ above.    The impact of these
choices is not generally as significant as is the choice of the threshold $E$.
\end{remark}

\subsection{Results of the feed-back based projection}
\label{sec:resfeed}
We applied the algorithm, with the threshold set at $E=0$,
to a random sample of 100,000 curves satisfying Angle Constraints 1 and 2.  We implemented the algorithm by computing  $e_\lambda$ at $t=0, 0.1, 0.2, \dots, 1$, taking the solution with lowest $e_\lambda$, and then repeating the procedure once. For the projected sample we obtained for $e_\lambda$:
\[
\hbox{mean}= 0.048,  \quad
\hbox{median} = 0.042, \quad
\hbox{max}  =  0.26.
\]
The 99th percentile is $e_\lambda =0.15$, and  $99.99\%$ of curves have $e_\lambda$ less than $0.22$.
Hence, one could set the target value $E$ to $0.22$ and expect the result to always have $e_\lambda \leq E$,
or set the target value at $0.15$ and expect this result $99\%$ of the time.

\begin{figure}[htb]
    \centering
    \begin{tabular}{cc}
												  \includegraphics[height=55mm]{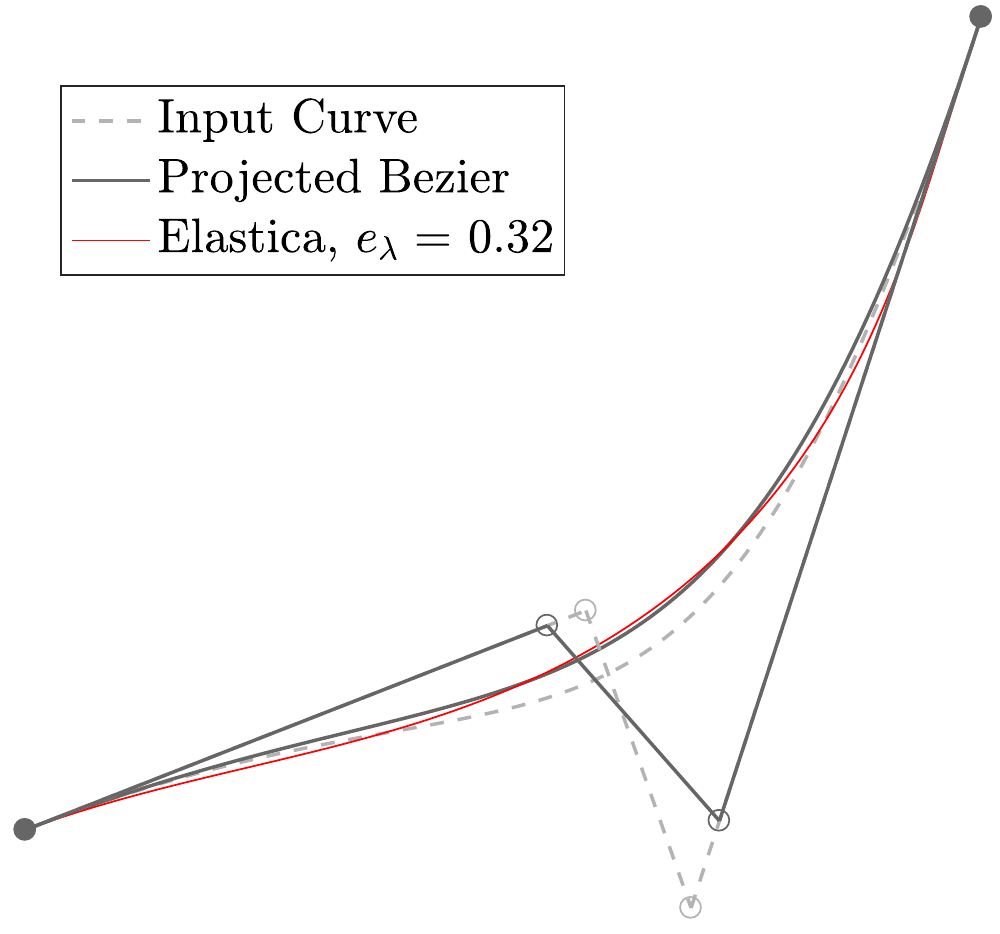}  \quad \quad &\quad \quad
						  \includegraphics[height=55mm]{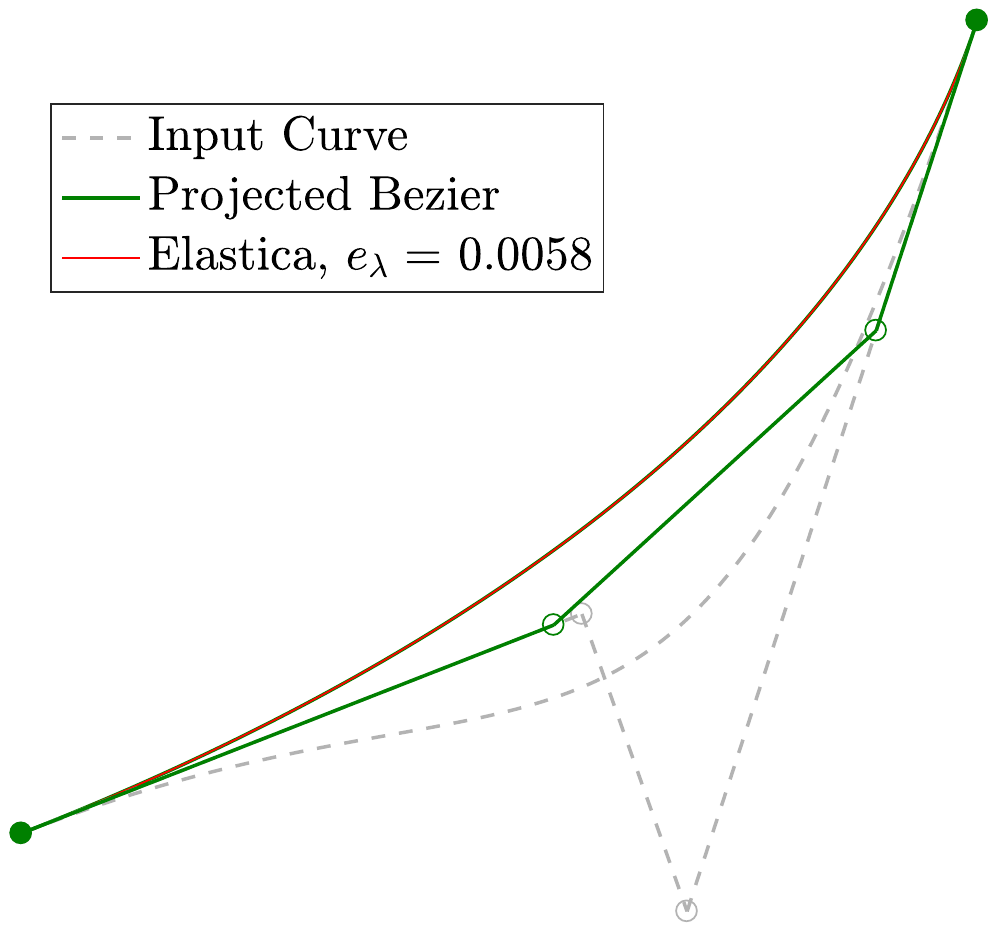}  
        \end{tabular}
    \caption{Result of the projection algorithm applied with two different
		choices of threshold level $E$ for $e_\lambda$.  This example changes from inflectional
		to non-inflectional.
		The red curve is the elastica obtained from
		the first guess.} 
    \label{fig:exampleprojectionsb}
\end{figure}

\begin{figure}[htb]
    \centering
				   \begin{tabular}{ccc}			
					  \includegraphics[height=37mm]{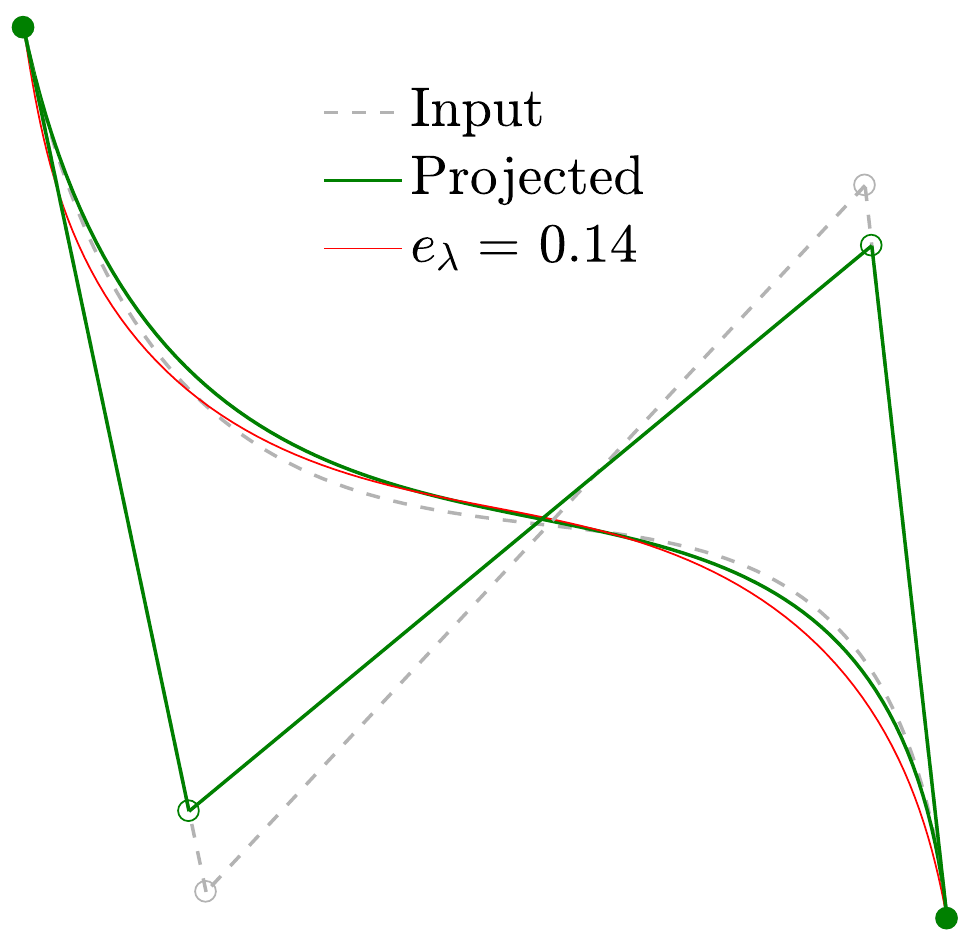}   
						   \quad & 
						  \includegraphics[height=37mm]{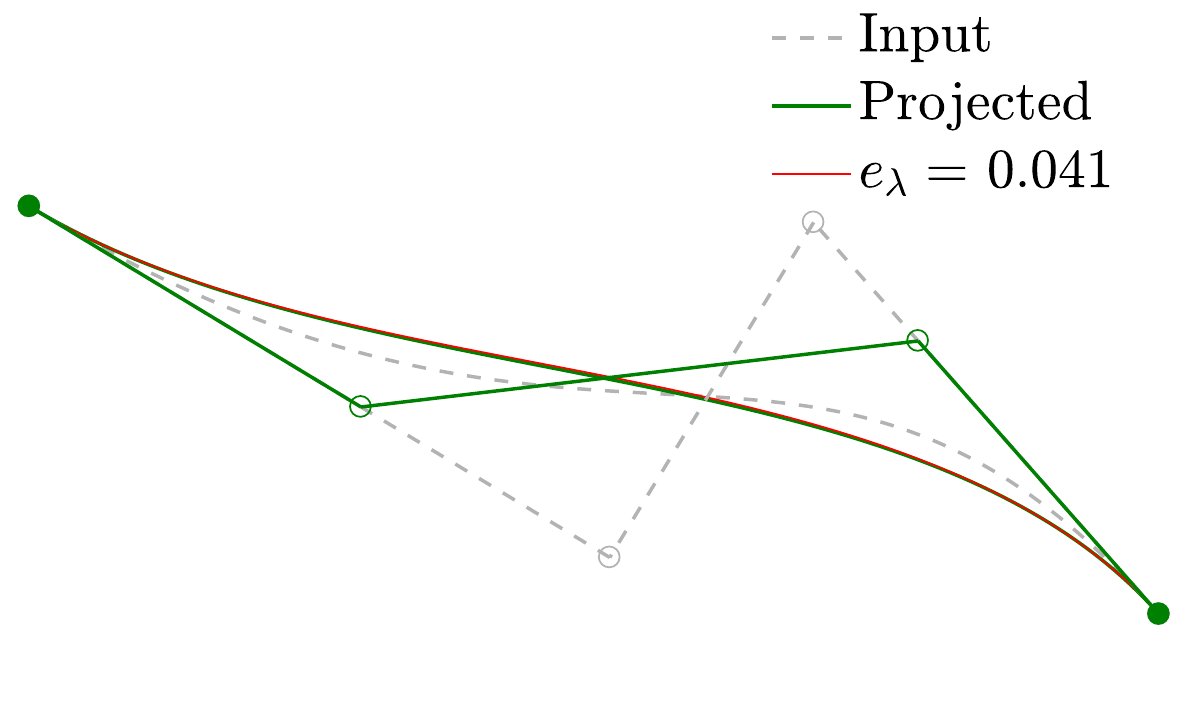}  & \quad
						  \includegraphics[height=37mm]{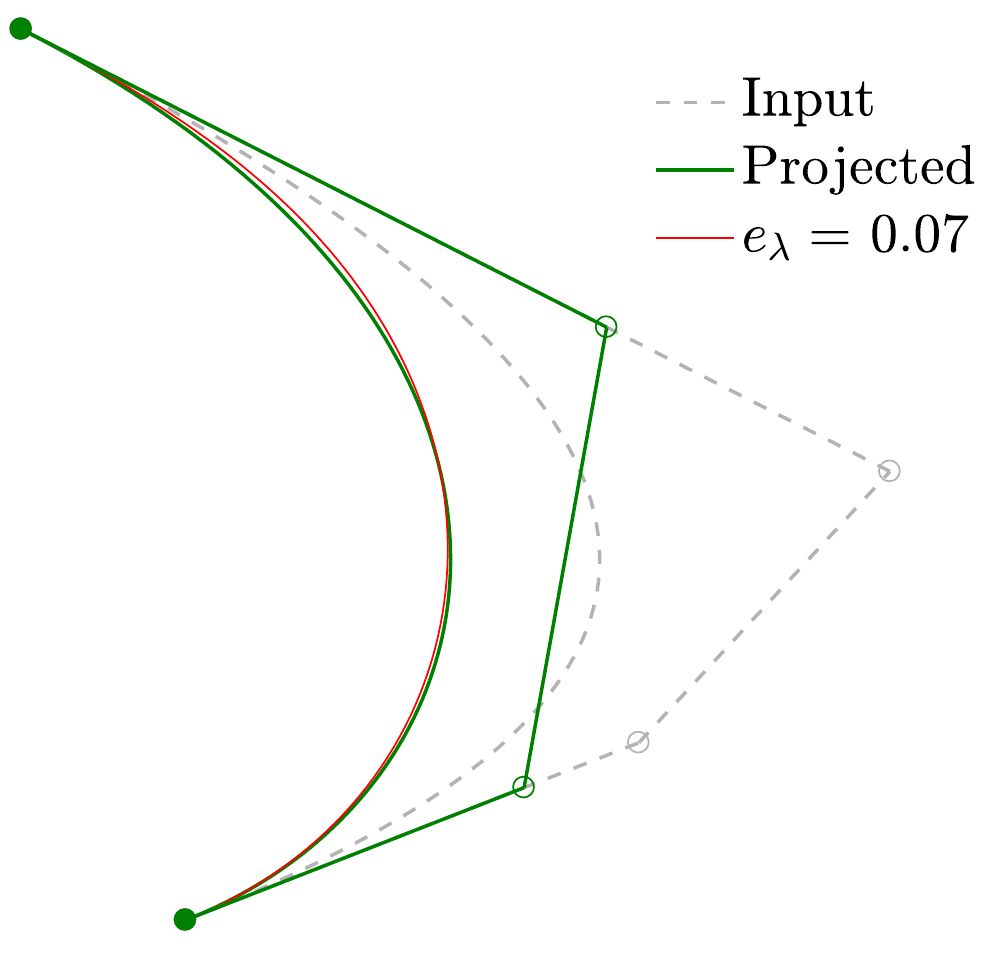} 

  \vspace{1ex} \\				
 \end{tabular}
    \caption{Examples of projections.  The elastica shown is the first guess.} 
    \label{fig:exampleprojections2}
\end{figure}

Representative examples satisfying Angle Constraints 1 and 2 are shown in Figures \ref{fig:exampleprojectionsb}
and \ref{fig:exampleprojections2}.  The elastic curve obtained just from the first guest algorithm in 
\cite{Brander}, (i.e., before any optimization is applied) is also plotted.  Note that, if a further optimization is
applied to these elastic curves, then they will approximate the projected B\'ezier curve significantly better.

\subsection{Extending the feedback projection to general B\'ezier curves}
A further advantage of the feedback method is that we can apply this projection to any B\'ezier curve,
regardless of the angle constraints, and use the feedback to decide whether or not the projected
B\'ezier curve is close to an elastica.  We took a sample of 83,500 arbitrary B\'ezier curves, each given by
4 randomly chosen control points in the unit disc.  Of these, $55\%$ satisfy both Angle Constraints 1 and 2.
Applying the projection algorithm (with target $E=0$) to the curves that do \emph{not} satisfy the angle constraints resulted in a mean for $e_\lambda$ of $0.4$ median of $0.3$, and maximum approximately 1. So the projection
should not be applied to B\'ezier curves with no angle constraints at all.

If we replace $\pi/3$ in Angle Constraint 1 by the value $\pi/4$, and replace the value $0.4\pi$ by
the value $0.6\pi$ in Angle Constraint 2, we still obtain good results of
\[
\hbox{mean}= 0.06,  \quad
\hbox{median} = 0.043, \quad
\hbox{max}  = 0.42,
\]
for $e_\lambda$ in the projected curve.  Extending further still to $\pi/6$ and $0.75\pi$
respectively produced a mean, median and maximum of $0.1$, $0.06$ and $0.71$ respectively,
which  is still a useable design space, provided that the feedback is used to reject the input curves
that exceed some desired threshold (e.g. those curves above 0.4). 

\begin{figure}[htb]
    \centering
    \begin{tabular}{cccc}
		\includegraphics[height=34mm]{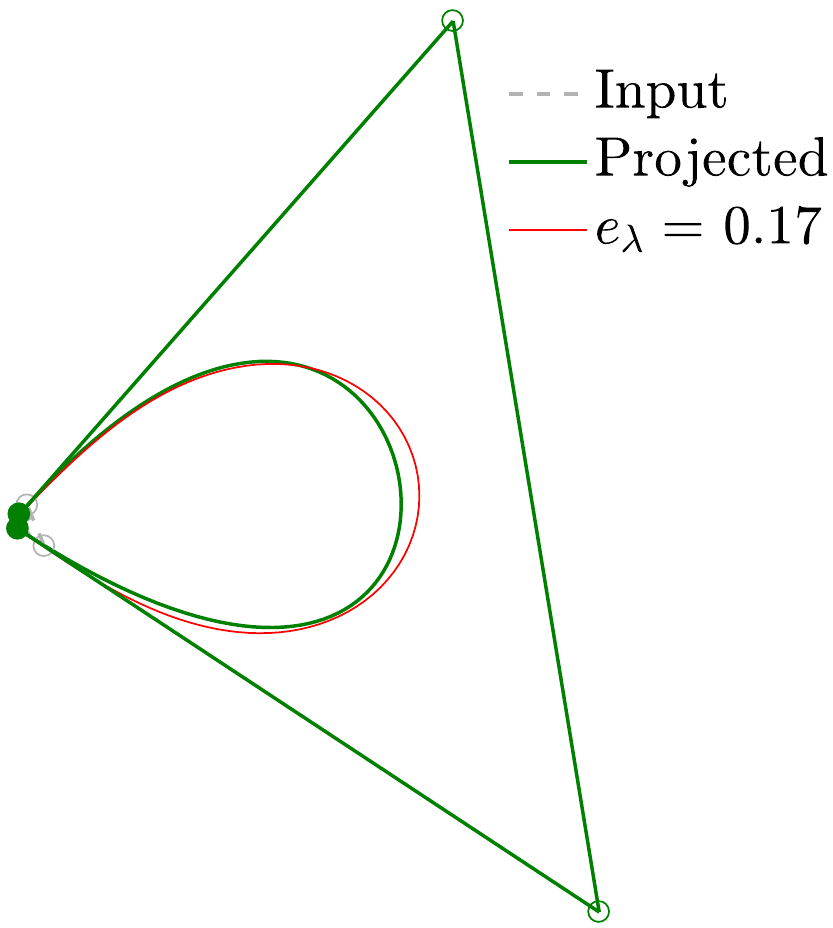}   \,\, & \,\,
			  \includegraphics[height=34mm]{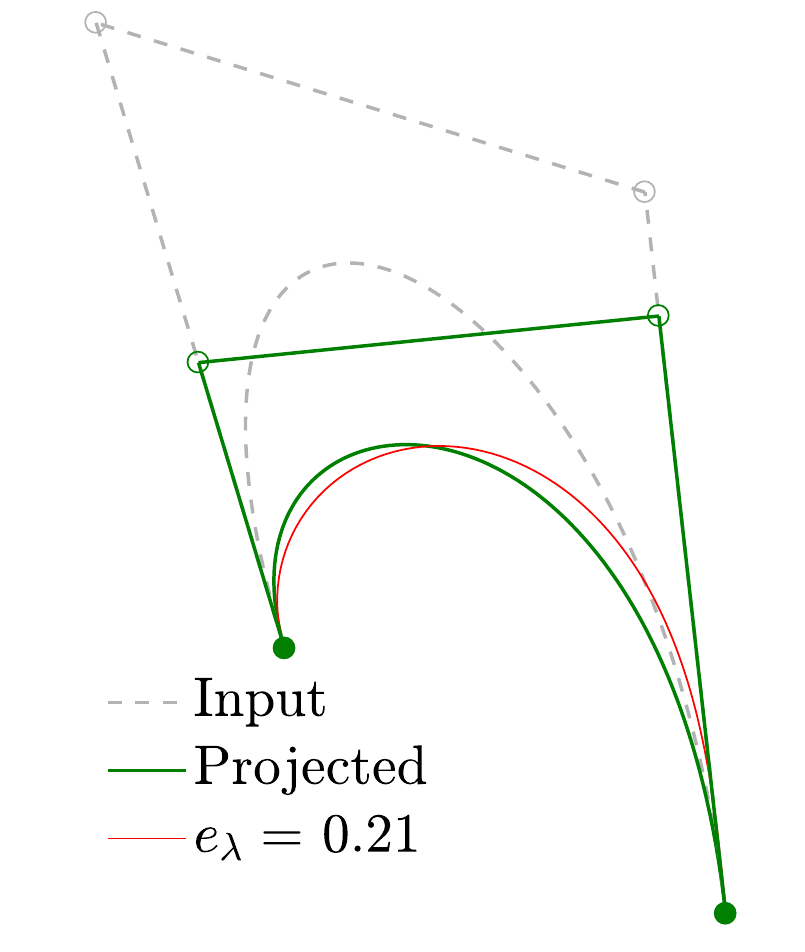} \,\, & \,\,
						\includegraphics[height=34mm]{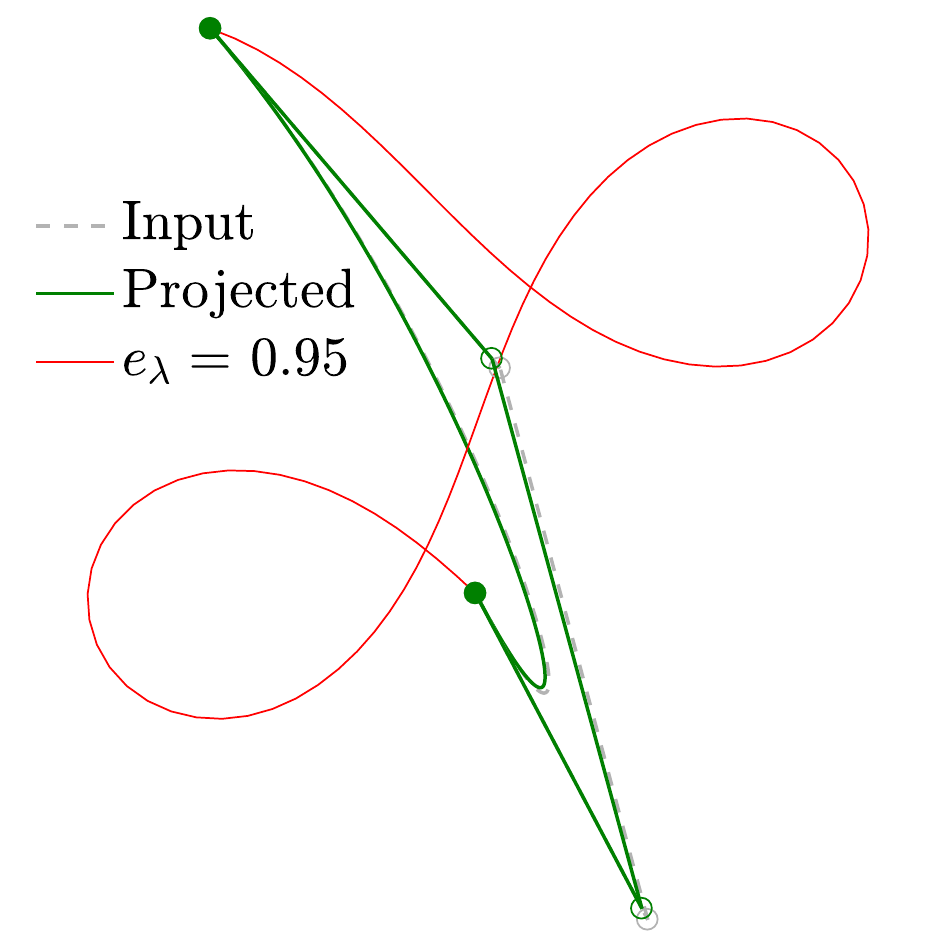}   \,\, & \,\,
								\includegraphics[height=34mm]{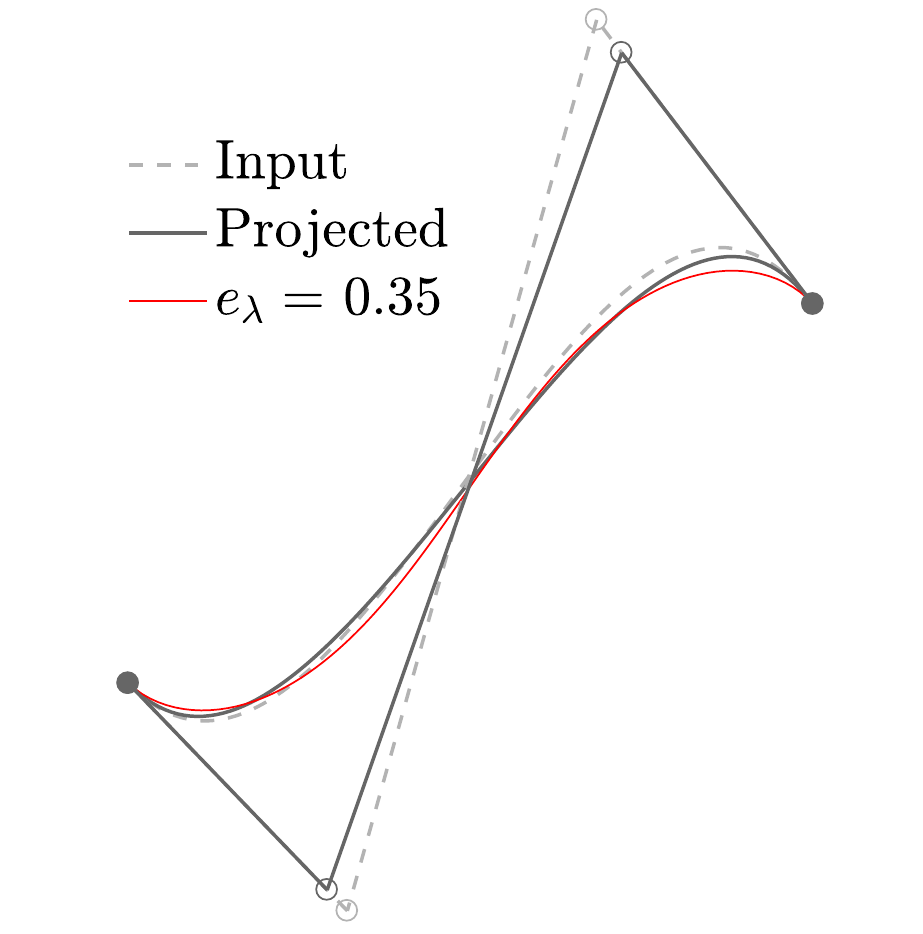}  
  \vspace{1ex} \\	
 \end{tabular}
    \caption{Examples that do not satisfy Angle Constraint 1 (left two), and Angle Constraint 2 (right two).} 
    \label{fig:exampleswide}
\end{figure}

Figure \ref{fig:exampleswide} shows some examples of the result of the projection applied to
curves that do not satisfy our original angle constraints.

\section{Conclusion and future work}
\label{sec:conclusion}
We have shown statistically that the $\lambda$-residual gives a convenient, easily computable, way
of measuring closeness to an elastic curve.   Based on this we have defined a reliable feedback-based
projection algorithm that takes an arbitrary cubic B\'ezier curve as input and adjusts the length to produce
a new cubic B\'ezier curve with the same endpoints and end-tangent angles.
  If the end-tangent angles of the input curve satisfy our angle constraints, then the
output curve is guaranteed to be close to an elastic curve.    For arbitrary B\'ezier curves, the end-tangent
angles can be rotated first if necessary.   An implementation in MATLAB of the feedback based projection
algorithm can be downloaded (at time of writing) from
\href{http://geometry.compute.dtu.dk/software/}{http://geometry.compute.dtu.dk/software/}
and tested.

A version of the projection algorithm described here has been 
incorporated into a hot-blade cutting design tool implemented in Rhino/Grashopper,
extending the work described in \cite{aag2016}.  This allows architects to design fabrication-ready
surfaces for hot-blade produced concrete casting.  The benefit of using elastica-like B\'ezier curves is that
the designed surface is visually the same as the  surface that will actually be produced.

\begin{figure}[htb]
    \centering
    \begin{tabular}{ccc}
			   \includegraphics[height=35mm]{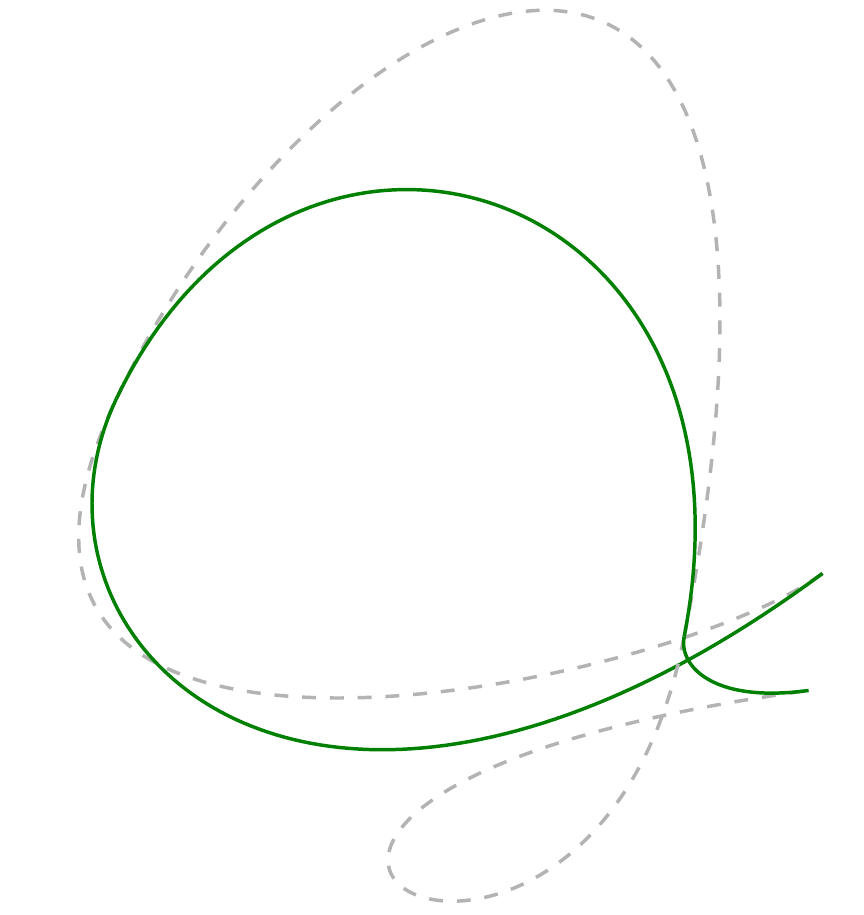}  \quad & \quad
								    \includegraphics[height=35mm]{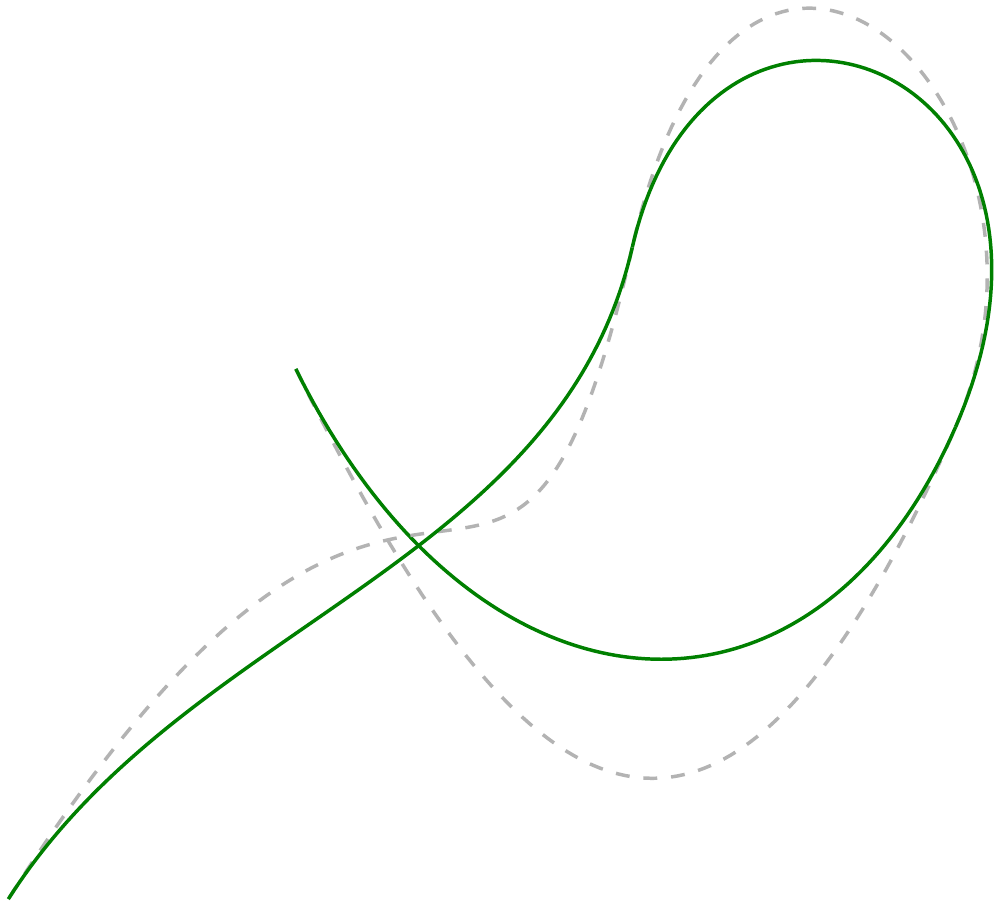}  \quad & \quad
				   \includegraphics[height=35mm]{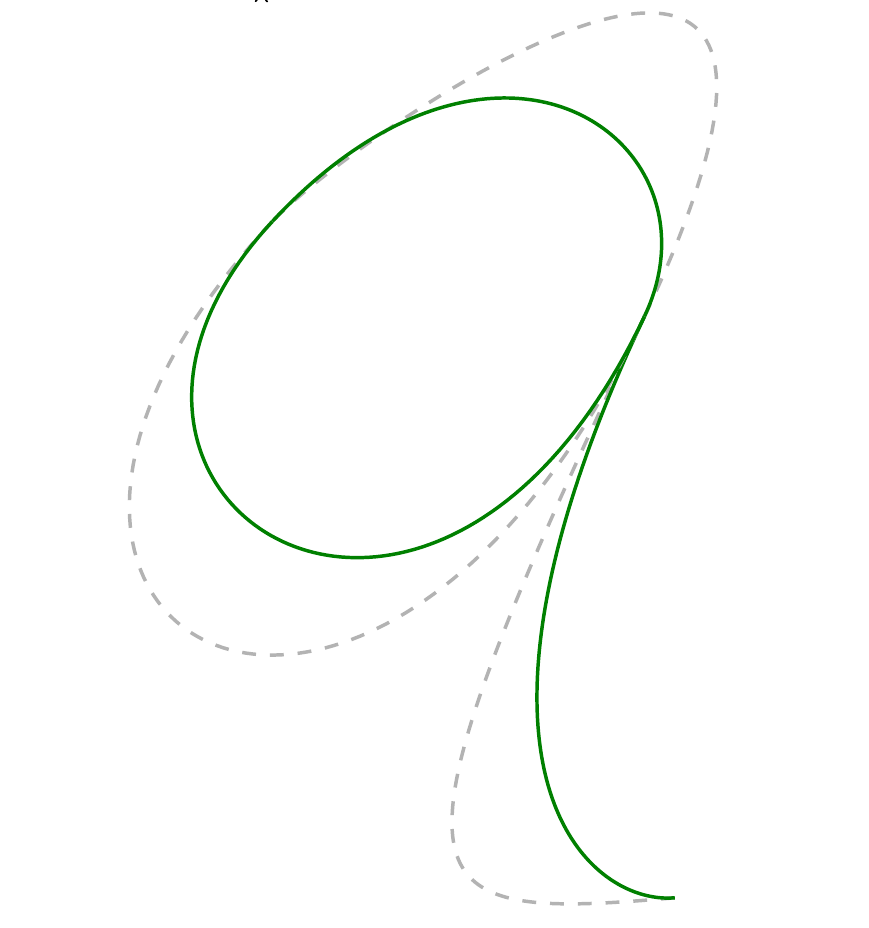}  
  \vspace{1ex} \\
 \end{tabular}
    \caption{Examples of $C^1$ cubic splines where each of 3 segments is close to an elastic curve.} 
    \label{fig:examplesplines}
\end{figure}

The projection algorithm can be used to construct $C^1$ splines that are close to elastic splines,
see  Figure~\ref{fig:examplesplines}.
However, for splines, a different method that allows end-points and end-tangents to be adjusted 
in addition to length is likely to be more suitable, and this will be studied in future work.

\section*{Acknowledgements}
Research partially supported by Innovation Fund Denmark, project number 91-2014-3.

\end{document}